\documentclass[12pt]{amsart}
\usepackage{latexsym,color,amsmath,amsthm,amssymb,amscd,amsfonts,enumitem}
\usepackage[bottom,symbol]{footmisc}
\usepackage{tikz}
\usepackage[normalem]{ulem}

\newtheorem{thm}{Theorem}[section]
\newtheorem*{thm*}{Theorem}

\newtheorem{cor}[thm]{Corollary}
\newtheorem{lem}[thm]{Lemma}
\newtheorem{prop}[thm]{Proposition}
\newtheorem{ques}[thm]{Question}
\newtheorem{conj}[thm]{Conjecture}
\newtheorem*{conj*}{Conjecture}

\theoremstyle{definition}
\newtheorem{rem}[thm]{Remark}

\newtheorem*{definition*}{Definition}

\newtheorem*{rems*}{Remarks}

\newcommand{\R}{{\mathbb{R}}}
\newcommand{\C}{{\mathbb{C}}}
\newcommand{\Z}{\mathbb{Z}}

\title {From Lagrangian Products to Toric Domains via the Toda Lattice}
	\author[Y. Ostrover]{Yaron Ostrover}
 \address{School of Mathematical Sciences, Tel Aviv University, Israel.}
 \email{ostrover@tauex.tau.ac.il}
 \author[V. G. B. Ramos]{Vinicius G. B. Ramos} 
  \address{Institute for Advanced Study, USA and Instituto de Matem\'atica Pura e Aplicada, Brazil.}
 \email{vgbramos@impa.br}
 
 \author[D. Sepe]{Daniele Sepe}
\address{Escuela de Matem\'aticas, Universidad Nacional de Colombia sede
  Medell\'in, Colombia, and Instituto de Matem\'atica e Estat\'istica, Departamento de
  Matem\'atica Aplicada (GMA), Universidade Federal Fluminense, Brazil.}
\email{dsepe@unal.edu.co}

\begin{document}
	
	\maketitle
	\begin{abstract}  
In this paper we use the periodic Toda lattice to show that certain Lagrangian product configurations in the classical phase space are symplectically equivalent to toric domains.
In particular, we prove that the Lagrangian product of a certain simplex and the Voronoi cell of the root lattice $A_n$ is symplectically equivalent to a Euclidean ball. As a consequence, we deduce that the Lagrangian product of an equilateral triangle and a regular hexagon is symplectomorphic to a Euclidean ball in dimension 4.
	\end{abstract}
	
\section{Introduction and main results} \label{sec-int}

The question of what can be said about the symplectic image of the Euclidean ball has its roots in Gromov's celebrated non-squeezing theorem~\cite{Gr},  
one of the most striking and fundamental results in symplectic topology. 
In this note we use the periodic Toda lattice to study this question for some Lagrangian products in the classical phase space. More precisely, given $A,B \subset {\mathbb R}^n$ with non-empty interior, the Lagrangian product of $A$ and $B$ is 
$$ A \times_L B := \{(x_1,y_1,\ldots,x_n,y_n) \in {\mathbb R}^{2n} \ | \ (x_1,\ldots,x_n) \in A \ {\rm and} \ (y_1,\ldots,y_n) \in B\},
$$
endowed with the restriction of the standard symplectic form $\omega_0 = \sum_{i=1}^n dx_i \wedge dy_i$. 

\begin{ques} \label{ques-Lag-prod-symp-to-B} 
For which $A, B \subset \R^n$ is $A \times_L B$ a symplectic ball?~\footnote{This question has two, \textit{a priori} distinct, interpretations: either as a symplectic equivalence (see Theorem \ref{thm-simplex-RD}), or as a symplectomorphism of the interiors.}
\end{ques}
 Lagrangian products arise naturally in the study of the dynamics of Minkowski billiards \cite{AAO}, and play an important role in the connection between Viterbo's volume-capacity conjecture~\cite{V} and Mahler’s conjecture~\cite{Ma} established in~\cite{AAKO}.
This motivated the suggestion to study symplectic embedding problems for Lagrangian products (see \cite[Section 5]{O}), which started in \cite{R}, and was subsequently continued in \cite{Os_Ra,RS}. 

To the best of our knowledge, up to scaling and affine symplectic transformations, there are not too many cases for which the answer to Question~\ref{ques-Lag-prod-symp-to-B} is known. A well-known example is the Lagrangian product $\Diamond^n \times_L \Box^n$, where $\Diamond^n$ is the 
standard octahedron in ${\mathbb R}^n$ (i.e., the convex hull of $\{\pm e_1, \pm e_2, \ldots, \pm e_n\}$), and $\Box^n = [0,1]^n$ is the standard cube in $\R^n$. The symplectomorphism between the interiors for $n=2$ appeared in \cite{LMS} (the proof therein uses a deep 4-dimensional result of McDuff~\cite{Mc}), and more recently in \cite{RS} for any dimension (cf. \cite{Schlenk}).
For the Lagrangian product $\triangle^n \times_L \Box^n$, where $\triangle^n$ is the standard simplex in ${\mathbb R}^n$ (i.e.,  
the convex hull of $\{0, e_1, e_2, \ldots, e_n\}$), the following is known: $\triangle^n \times_L \Box^n$ embeds symplectically into the Euclidean ball of an arbitrarily close volume, and vice versa (see \cite{LMS,Sc,Tr}, cf.~\cite{Ru}). Moreover, if $n=2$, the results in \cite{Mc} implies that the interior of $\triangle^2 \times_L \Box^2$ is symplectomorphic to an open ball. We observe that one can also prove the existence of a symplectomorphism\footnote{Schlenk communicated to us that, upon choosing the above symplectic embeddings carefully, \cite[Theorem 4.3]{PVN} yields the desired symplectomorphism.} between the interiors of $\triangle^n \times_L \Box^n$ and a Euclidean ball  using the strategy of \cite{RS} with a suitable integrable system. 


In this paper we provide another  
answer to Question~\ref{ques-Lag-prod-symp-to-B}. We use the periodic Toda lattice to show that the Lagrangian product of a certain simplex 
and the Voronoi cell of the root lattice $A_n$ (see \cite[Chapter 21, Section 3.B]{CSbook}) 
is symplectically equivalent to a Euclidean ball. To state the result precisely, in what follows we identify $\R^{n-1}$ with the hyperplane 
$$W^{n-1}=\{(x_1,\dots,x_{n})\in\R^{n}\mid x_1+\dots+x_{n}=0\} \subset \R^{n}.$$
Let $S_{n}$ be the symmetric group on $n$. The above hyperplane is invariant under the $S_n$-action on $\R^{n}$ that permutes the coordinates. We consider the open sets 
\begin{equation}
\label{eq:100}
   \begin{aligned} 
\mathfrak{S}^{n-1} & = \left\{ (x_1,\ldots,x_{n}) \in W^{n-1} \ | \  x_{i} - x_{i+1} < 1 \ {\rm for} \ 1\leq i \leq n, \,  x_{n+1} := x_1 \right\} ,\\
\mathfrak{P}^{n-1}&= \left\{ (x_1,\ldots,x_{n}) \in W^{n-1} \ | \  \max(x_1,\ldots,x_{n})-\min(x_1,\ldots,x_{n})< 1 \right\}.\end{aligned} 
\end{equation}
The set  $\mathfrak{S}^{n-1}$ is the interior of the $(n-1)$-simplex given by the convex hull of the $n$ cyclic permutations of the vector 
$$ \Bigl (-{\frac {n-1} 2}, -{\frac {n-1} 2}+1, \ldots, {\frac {n-1} 2}-1,{\frac {n-1} 2} \Bigr ) \in W^{n-1}.$$
Note that the barycenter of $\mathfrak{S}^{n-1}$ is at the origin. The set $\mathfrak{P}^{n-1}$ is the convex hull of the $2n$ cyclic permutations of the vectors
\[\left(\frac{1}{n},\dots,\frac{1}{n},-\frac{n-1}{n}\right)\text{ and }\left(-\frac{1}{n},\dots,-\frac{1}{n},\frac{n-1}{n}\right).\]
The polytope $\mathfrak{P}^{n-1}$ can be described in terms of a difference body as we explain below. Consider the simplex 
$$\triangle_W^{n-1} = {\rm Conv} \Bigl \{ e_1, \ldots , e_{n-1}, - \sum_{i=1}^{n-1} e_i \Bigr \} \subset W^{n-1},$$
\noindent
where $\mathrm{Conv}$ denotes the convex hull.
\vspace{0.1cm}
\noindent
Note that $\triangle_W^{n-1}$ also has its barycenter at the origin. By definition of $\triangle_W^{n-1}$, the support function of the difference body $ \triangle_W^{n-1}-\triangle_W^{n-1}$ is  given by
$$h_{\triangle_W^{n-1}-\triangle_W^{n-1}}(x) =  \max_{1 \leq i \leq n-1} \Bigl \{x_i, -\sum_{i=1}^{n-1} x_i \Bigr \} - \min_{1 \leq i \leq n-1} \Bigl \{x_i, -\sum_{i=1}^{n-1} x_i \Bigr \},$$
where we recall that the support function of a convex body $ K \subset {\mathbb R}^{n-1}$ is  given by $h_K(x) = \sup \{ \langle x,y \rangle <1 \mid y\in K \}$. Hence, $\mathfrak{P}^{n-1}$ is the polar body of $ \triangle_W^{n-1}-\triangle_W^{n-1}$. Such polytopes are sometimes known as ``generalized rhombic dodecahedra" (see, e.g.,~\cite{Z}). In particular, for $n=3$, the set $\mathfrak{P}^2$ is a regular hegaxon and for $n=4$, the set $\mathfrak{P}^3$ is a rhombic dodecahedron (see Figure~\ref{fig:RD} below). Since the two simplices $\mathfrak{S}^{n-1}$ and $\triangle_W^{n-1}$ share the same barycenter, one can construct a map $A \in {\rm SL}(n)$ such that  $\mathfrak{P}^{n-1} $ is the dual of the difference body of $A \mathfrak{S}^{n-1}$.

\begin{figure}[h]
    \centering
    \includegraphics[scale=0.2]{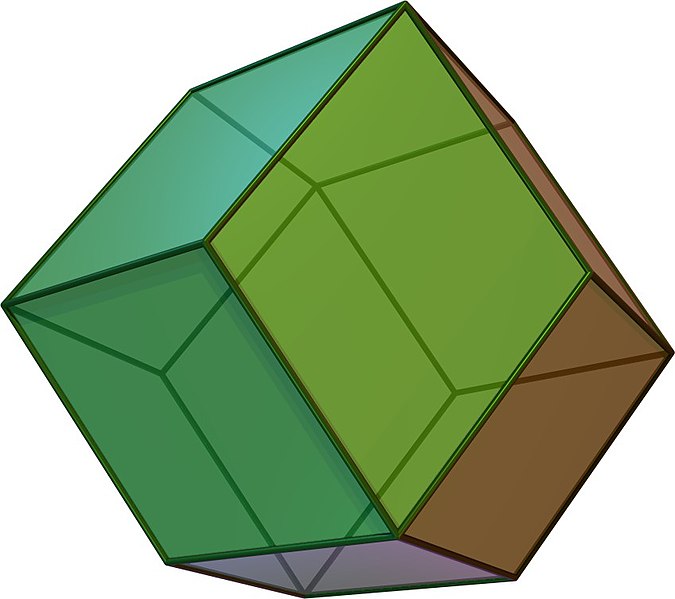}
  \caption{A rhombic dodecahedron (the dual of the difference body of a  simplex)\protect\footnotemark.}
    \label{fig:RD}
\end{figure}

\noindent

The subspace  $V^{2n-2}=W^{n-1}\times_L W^{n-1}\subset\R^{2n}$ is symplectic. Hence, in what follows we may identify $V^{2n-2}$ with $\R^{2n-2}$ symplectically, so that we may consider Lagrangian products of the form $A \times_L B$, where $A, B \subset W^{n-1}$. Moreover, we use the following notation for a Euclidean ball in $\R^{2n-2}$: 
$$B^{2n-2}(a) := \{(x_1,y_1,\dots,x_{n-1},y_{n-1})\in\R^{2n-2}\mid \pi\left(x_1^2+y_1^2+\dots+x_{n-1}^2+y_{n-1}^2\right)<a\}.$$
The precise (partial) answer to Question \ref{ques-Lag-prod-symp-to-B} is the following result.

    \begin{thm} \label{thm-simplex-RD}
The Lagrangian product $\mathfrak{S}^{n-1} \times_L \mathfrak{P}^{n-1}$ is symplectically equivalent to $B^{2n-2}(n)$, i.e., for every $\epsilon>0$, there exist symplectic embeddings
\[(1-\epsilon)\left(\mathfrak{S}^{n-1} \times_L \mathfrak{P}^{n-1} \right)\hookrightarrow B^{2n-2}(n)\hookrightarrow(1+\epsilon)\left(\mathfrak{S}^{n-1} \times_L \mathfrak{P}^{n-1}\right).\]
\end{thm}
 \footnotetext{The figure by w:en:User:Cyp@wikimedia, distributed under CC BY-SA 3.0 license.}

We observe that, for $n=3$, Theorem~\ref{thm-simplex-RD} also implies a partial answer to  \cite[Question 4.4]{Bal} (cf. \cite[Theorems 1.2 and 1.3]{Ru}). Indeed, in this case, it follows from this theorem that the Lagrangian product of an equilateral triangle $T$ 
and the hexagon $T-T$ rotated by $\pi/6$ degrees is symplectomorphic to a Euclidean ball in ${\mathbb R}^4$. More precisely, the following holds.
\begin{figure}[ht] 
\begin{center}
  \begin{tikzpicture}[scale=1]

\path coordinate (q1) at (0,1.414) coordinate (q2) at (-1.224,-0.707) coordinate (q3) at  (1.224,-0.707);

\draw[line width=0.35mm] (q1) -- (q2);
\draw[line width=0.35mm] (q2) -- (q3);
\draw[line width=0.35mm] (q3) -- (q1);

\draw[thick,->,blue] (-2,0)--(2,0) node[below] {$x_1$}; 
    \draw[thick,->,blue] (0,-1.5)--(0,2) node[left] {$x_2$}; 
    \draw (0,-1.6) node[below] {$\widetilde {\mathfrak{S}}^2$};

\filldraw [black]
 (q1) circle (1pt) node[right] {{\footnotesize $(0,\sqrt{2})$}}
  (q2) circle (1pt) node[left] {{\footnotesize $(-{\frac {\sqrt 3} {\sqrt 2}},{\frac {-1} { \sqrt{2}}})$}}
      (q3) circle (1pt) node[right] {{\footnotesize $({\frac {\sqrt 3} {\sqrt 2} },{\frac {-1} { \sqrt{2}}})$}};

      \begin{scope}[xshift=7cm]

\path coordinate (p1) at (5*0.272,0) coordinate (p2) at (5*0.136,5*0.235) coordinate (p3) at (-0.136*5,0.235*5) coordinate (p4) at (-0.272*5,0)
coordinate (p5) at (-0.136*5,-0.235*5) coordinate (p6) at  (0.136*5,-0.235*5);
\draw[line width=0.35mm]  (p1) -- (p2);
\draw[line width=0.35mm]  (p2) -- (p3);
\draw[line width=0.35mm]  (p3) -- (p4);
\draw[line width=0.35mm]  (p4) -- (p5);
\draw[line width=0.35mm]  (p5) -- (p6);
\draw[line width=0.35mm] (p6) -- (p1);

\filldraw [black]

(p1) circle (1pt) node[below right] {{\footnotesize $({\frac {\sqrt{6}} {3} },0)$}}

      (p2) circle (1pt) node[right] {{\footnotesize $({\frac 1 { \sqrt{6}}},{\frac {1} {\sqrt{2}}})$}};

\draw[thick,->,blue] (-2,0)--(2,0) node[right] {$y_1$}; 
    \draw[thick,->,blue] (0,-1.5)--(0,2) node[left] {$y_2$}; 
\draw (0,-1.6) node[below] {$\widetilde {\mathfrak{P}}^2$};

\end{scope}
\end{tikzpicture}

\caption{The product $\widetilde{\mathfrak{S}}^2  \times_L \widetilde{ \mathfrak{P}}^2$ is symplectomorphic to a Euclidean ball.}\label{fig:test}
\end{center}
\end{figure}
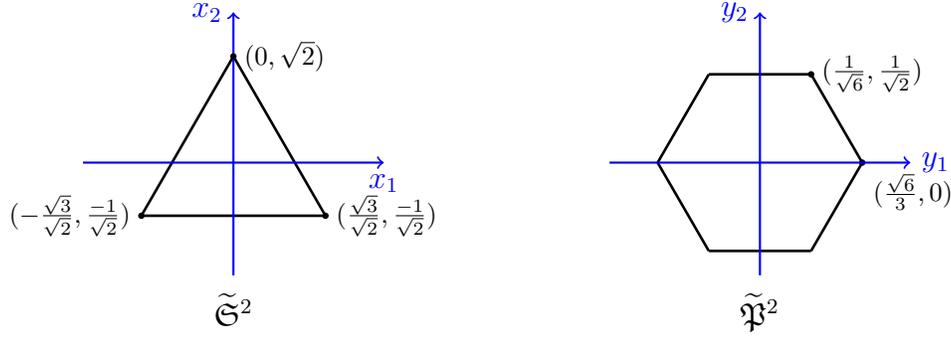
 \begin{cor}\label{cor}
 The interior of the Lagrangian product of the equilateral triangle and the regular hexagon of Figure \ref{fig:test} is symplectomorphic to the interior of a symplectic ball.
 \end{cor}
\begin{proof}
If $n=2$, then $\mathfrak{S}^2$ is the interior of the equilateral triangle  
$$ {\rm Conv} \bigl \{ (-1,0,1),(1,-1,0),(0,1,-1) \bigr \} \subset W^2,$$
and $\mathfrak{P}^2 $ is the interior of the regular hexagon  
$$ \mathfrak{P}^2 = {\rm Conv} \bigl \{ ({\tfrac {-2} 3},{\tfrac 1 3},{\tfrac 1 3}),({\tfrac 1 3},{\tfrac {-2} 3},{\tfrac 1 3}),({\tfrac 1 3},{\tfrac 1 3},{\tfrac {-2} 3}) , ({\tfrac {2} 3},{\tfrac {-1} 3},{\tfrac {-1} 3}),({\tfrac {-1} 3},{\tfrac {2} 3},{\tfrac {-1} 3}),({\tfrac {-1} 3},{\tfrac {-1} 3},{\tfrac {2} 3})\bigr \} \subset W^2.$$ 
Note that the Lagrangian product $\mathfrak{S}^2 \times_L \mathfrak{P}^2$ is symplectomorphic to $\widetilde {\mathfrak{S}}^2 \times_L \widetilde{ \mathfrak{P}}^2$, where $\widetilde {\mathfrak{S}}^2 = A \mathfrak{S}^2$, $\widetilde {\mathfrak{P}}^2 = (A^T)^{-1} \mathfrak{P}^2$,
and $A \in {\rm O}(3)$ is the orthogonal transformation that rotates the plane $x_1+x_2+x_3=0 $ in ${\mathbb R}^3$ to the plane $x_3=0$ so that $\widetilde{\mathfrak{S}}^2$ and $\widetilde{\mathfrak{P}}^2$ are as in Figure \ref{fig:test}. Hence, by Theorem \ref{thm-simplex-RD}, the Lagrangian product  $\widetilde{ \mathfrak{S}}^2  \times_L \widetilde {\mathfrak{P}}^2$ is symplectically equivalent to a Euclidean ball. By \cite[Lemma 4.3]{LMS}, the result follows.
\end{proof}
 
Theorem~\ref{thm-simplex-RD} follows from the main result of this paper, Theorem \ref{thm-lag-prod-as-toric-domain}. Intuitively speaking, the latter states that the Lagrangian product of $\mathfrak{S}^n$ and a suitable $S_{n}$-invariant region in $W^{n-1}$
is symplectically equivalent to a toric domain, i.e.,  a subset of $\R^{2n-2}$ of the form 
\[X_\Omega=\left\{(x_1,y_1,\dots,x_{n-1},y_{n-1})\in\R^{2n-2}\mid \left(\pi(x_1^2+y_1^2),\dots,\pi(x_{n-1}^2+y_{n-1}^2)\right)\in\Omega\right\},\]
where $\Omega\subset \R^{n-1}_{\ge 0}$ is a relatively open set. The terminology comes from the fact that, given $\Omega \subset \R^{n-1}_{\ge 0}$, the subset $X_{\Omega} \subset \R^{2n-2}$ is invariant under the standard $(S^1)^{n-1}$-action on $\R^{2n-2} \simeq \C^{n-1}$, where $S^1 = \R/\Z$. A Euclidean ball in $\R^{2n-2}$ is an example of a toric domain: given $a > 0$, $B^{2n-2}(a) = X_{T(a)}$, where $T(a)=\{(r_1,\dots,r_{n-1})\in\R^{n-1}_{\ge 0}\mid r_1+\dots+r_{n-1}<a\}$. 

In order to state Theorem \ref{thm-lag-prod-as-toric-domain}, we need to introduce some more notions.
In what follows we restrict our discussion to subsets of $W^{n-1}$ that are $S_{n}$-invariant. To any such subset, we may naturally associate a subset of $\R^{n-1}_{\ge 0}$ as follows: For any $(y_1,\ldots, y_{n}) \in \R^{n}$, there exists $\sigma \in S_{n}$ such that $y_{\sigma(1)}\le \dots \le y_{\sigma(n)}$. Although such a permutation may not be unique, the vector 
\begin{equation} \label{eqn-def-rho}\rho(y_1,\dots,y_{n}):=n \bigl (y_{\sigma(2)}-y_{\sigma(1)},y_{\sigma(3)}-y_{\sigma(2)},\dots,y_{\sigma(n)}-y_{\sigma(n-1)} \bigr) \end{equation}
does not depend on the choice of permutation. Hence, the map $\rho : \R^{n} \to \R^{n-1}_{\geq 0}$ is well-defined. To an $S_{n}$-invariant subset $A \subset W^{n-1}$ we associate $\rho(A) \subset \R^{n-1}_{\geq 0}$. We remark that if $A$ is $S_{n}$-invariant, the set $\rho(A)$ determines $A$ uniquely. 
Finally, we also say that a set $A\subset W^{n-1}$ is star-shaped if $rA\subset A$ for all $0\le r \le 1$. The main result of our paper is the following.

\begin{thm} \label{thm-lag-prod-as-toric-domain}
Let $A\subset W^{n-1}$ be a subset that is $S_{n}$-invariant, open, bounded and star-shaped. Then the Lagrangian product $\mathfrak S^{n-1} \times_L A$ is symplectically equivalent to the toric domain $X_{\rho(A)}$.

\end{thm}
Theorem \ref{thm-simplex-RD} is a simple consequence of Theorem \ref{thm-lag-prod-as-toric-domain}. 

\begin{proof}[{\bf Proof of Theorem \ref{thm-simplex-RD}}]
The subset $\mathfrak{P}^{n-1} \subset W^{n-1}$ of \eqref{eq:100} is $S_{n}$-symmetric, relatively open, bounded and convex. 
For $(y_1,\dots,y_{n})\in \mathfrak{P}^{n-1}$,  set $(r_1,\dots,r_{n-1}) :=\rho(y_1,\dots,y_{n})$. By the definition of $\rho$, 
\[r_1+\dots+r_{n-1}=n\left(\max\{y_1,\dots,y_{n}\}-\min\{y_1,\dots,y_{n}\}\right)<n.\]
Conversely, let $(r_1,\dots,r_{n-1})\in \R^{n-1}_{\ge 0}$ be such that $\sum_{i=1}^{n-1} r_i < n$. 
For $i=1,\dots,n$, we set \[y_i:=\frac{1}{n^2}\sum_{j=1}^{n-1} j r_j-\frac{1}{n}\sum_{j=i}^{n-1} r_j.\]
By a simple calculation, we have that $y_1+\dots+y_{n} = 0$, and
\[\begin{aligned} 
\max\{y_1,\dots,y_{n}\}-\min\{y_1,\dots,y_{n}\}&=y_{n}-y_1=\frac{1}{n}\sum_{j=1}^{n-1} r_j<1.\end{aligned}\]
Hence,
 $ \rho(\mathfrak{P}^{n-1})=T(n)$ and then $X_{\rho(\mathfrak{P}^{n-1})}=B^{2n-2}(n)$.
By Theorem \ref{thm-lag-prod-as-toric-domain}, the Lagrangian product 
$\mathfrak{S}^{n-1}\times_L \mathfrak{P}^{n-1}$ is symplectically equivalent to $B^{2n-2}(n)$. 
\end{proof}

While the strategy of the proof of Theorem \ref{thm-lag-prod-as-toric-domain} is somewhat similar to that of the main result in \cite{RS}, the ingredients are substantially more involved. The 1-parameter family of integrable systems that we use to construct the desired symplectic embeddings comes from the periodic Toda lattice (see \cite{avmm,Fl1,Fl2,Hen,To1} and references therein). It is known that this system admits global action-angle variables (see \cite{FM, HeKa,vM}). The proof of Theorem \ref{thm-lag-prod-as-toric-domain} hinges on showing certain properties of the global action variables as a certain parameter goes to infinity (see Proposition \ref{prop:lim}), and the combination of these with the results in \cite{FM, HeKa,vM}. In particular, we show that the standard billiard in (the closure of) $\mathfrak{S}^{n-1}$ can be thought of as a limit of a suitable 1-parameter deformation of the standard Toda lattice (see Lemma \ref{lemma:uc} and Remark \ref{rem:limit}). We believe that this result should be of independent interest.

\begin{figure}[h!]
    \centering
    \includegraphics[scale=0.6]{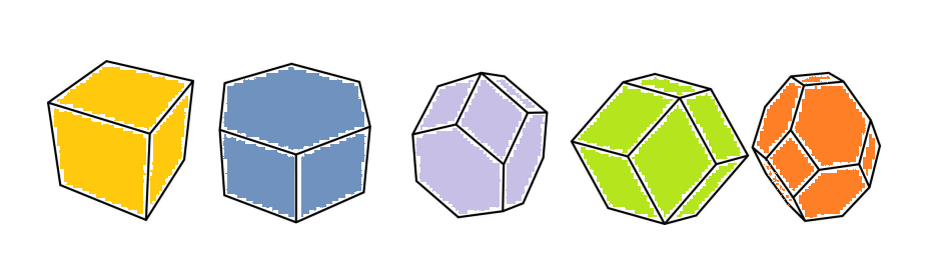}
    \caption{The five Fedorov solids (from left to right): the cube, the hexagonal prism, the elongated dodecahedron, the rhombic dodecahedron, and the truncated octahedron.}
    \label{fig:Fedorov}
\end{figure}

Since the proof of Theorem~\ref{thm-lag-prod-as-toric-domain} involves the periodic Toda Hamiltonian associated with the root system $A_n$, 
we speculate that the study of the more general periodic Toda potentials considered in \cite{Bog} and \cite[Chapter 9]{avmm} might lead to different solutions to Question~\ref{ques-Lag-prod-symp-to-B}. In fact, going even further, we expect that more interesting results in symplectic topology can be obtained by studying the relation between classical integrable systems and billiards. In particular, one possible candidate 
is the Lagrangian product of a simplex and a permutohedron in a certain position. Indeed, it was proved in~\cite{Bal}
that the Lagrangian product of a simplex and a permutohedron (properly related to each other) delivers equality in Viterbo's conjecture~\cite{V}. It is plausible that one can use the same strategy of proof as in Theorem~\ref{thm-lag-prod-as-toric-domain} to show that the above Lagrangian product is symplectically equivalent to a Euclidean ball, thus obtaining a complete answer to~\cite[Question 4.4]{Bal}. Moreover, in dimension three, motivated by Theorem~\ref{thm-simplex-RD}, we can generalize 
from~\cite[Question 4.4]{Bal}
as follows: 
recall that in dimension three there are five space filling polytopes that tile the space by translations (also known as parallelohedra or Fedorov solids, see \cite{CS} and Figure~\ref{fig:Fedorov}). 

\begin{conj} \label{conj-parallelohedra}
The Lagrangian product of a three dimensional simplex and any one of the five parallelohedra in an appropriate linear position is symplectically equivalent to a Euclidean ball. 
\end{conj}
Other than the two known cases of the cube and the rhombic dodecahedron, Conjecture~\ref{conj-parallelohedra} is supported by
the fact that all the aforementioned Lagrangian products deliver equality in Viterbo's conjecture~\cite{V}. Indeed, the case of the permutohedron follows from the main result of~\cite{Bal}, and in the remaining two cases the systolic ratio in Viterbo's conjecture can be explicitly computed based on the main result of~\cite{HK}.

\noindent{\bf Structure of the paper:} 
In Section~\ref{sec-examples} we show how to use Theorem~\ref{thm-lag-prod-as-toric-domain} to construct Lagrangian products that are symplectomorphic to ellipsoids and polydisks in dimension four (similar examples could be obtained also in higher dimensions). 
In Section~\ref{section-Toda} we recall some known facts regarding the periodic Toda lattice, and construct a certain symplectomorphism that plays 
an important role in the proof of  our main result. Section~\ref{section-from-Toda-to-billiards} is the technical heart of the paper and contains the proof of Theorem \ref{thm-lag-prod-as-toric-domain}.

\noindent{\bf Acknowledgements:} 
We are grateful to Alexey Balitskiy, Pazit Haim-Kislev and Felix Schlenk for useful discussions and valuable remarks. 
Y.O. was partially supported by the ISF grant No. 938/22. V.G.B.R. was partially supported by NSF grant DMS-1926686, FAPERJ grant JCNE E-26/201.399/2021 and a Serrapilheira Institute grant. D.S. was partially supported by FAPERJ grant JCNE E-26/202.913/2019 and by a CAPES/Alexander von Humboldt Fellowship for Experienced Researchers 88881.512955/2020-01. This study was financed in
part by the Coordenação de Aperfeiçoamento de Pessoal de Nível Superior -- Brazil
(CAPES) -- Finance code 001.

\section{Some examples in dimension four} \label{sec-examples}
In this section we use Theorem~\ref{thm-lag-prod-as-toric-domain} to construct Lagrangian products that are either symplectically equivalent or symplectomorphic to ellipsoids or polydisks in dimension four.

\subsection{Ellipsoids and Polydisks}

For $a,b>0$ consider the open ellipsoid and the polydisk in ${\mathbb C}^2$ given by   
$$ E(a,b) = \Bigl \{ (z_1,z_2) \, | \, {\frac {\pi|z_1|^2} {a} } + {\frac {\pi|z_2|^2} {b}} < 1 \Bigr \}, \ \ P(a,b) =  \Bigl \{ (z_1,z_2) \, | \, {\frac {\pi|z_1|^2} {a} } < 1, \, {\frac {\pi|z_2|^2} {b}} < 1 \Bigr \}.   $$
Both $E(a,b)$ and  $P(a,b)$ are toric domains: The moment map image of the former is the interior of the triangle with vertices $(0, 0), (a, 0)$, and $(0, b)$, while that of the latter is the interior of the rectangle with vertices $(0, 0), (a, 0), (0, b),(a, b)$. In analogy with Theorem~\ref{thm-simplex-RD}, in what follows we describe Lagrangian product configurations that are symplectically equivalent to either $E(a,b)$ or $P(a,b)$.

First, we deal with the case of the ellipsoid $E(a,b)$. Consider the set $$ \mathfrak{P}_E^2[a,b] :=  \bigsqcup_{\sigma \in S_3}  \sigma \cdot T[a,b] \subset W^2$$ that is the union of the six open triangles obtained by the action of $S_3$ on the interior of the triangle $$ T[a,b] := {\rm Conv} \bigl \{ (0,0,0),  ({\tfrac {-2a} 3},{\tfrac a 3},{\tfrac a 3}),({\tfrac {-b} 3},{\tfrac {-b} 3},{\tfrac {2b} 3}) \bigr \}.$$
Set $\widetilde{\mathfrak{S}}^2 = A \mathfrak{S}^2$, and $\widetilde{\mathfrak{P}}_E^2[a,b] = (A^{T})^{-1} \mathfrak{P}^2_E[a,b]$,
where $A \in {\rm O}(3)$ rotates the plane $x_1+x_2+x_3=0 $  to the plane $x_3=0$. 
Assume for simplicity that $0<a<b$. The closure of $\widetilde{\mathfrak{P}}_E^2[a,b]$ (in the $\{x_1,x_2\}$-plane) is the hexagon with a counterclockwise sequence of vertices given by
$$  \bigl \{ \bigl ({\tfrac {b\sqrt{6}} {3} },0 \bigr ), \bigl ({\tfrac {a} {\sqrt{6}}},{\tfrac {a} {\sqrt{2}} } \bigr ) , \bigl ({\tfrac {-b} {\sqrt{6}}},{\tfrac {b} {\sqrt{2}} } \bigr ) , \bigl (-{\tfrac {a\sqrt{6}} {3} },0 \bigr ),\bigl (-{\tfrac {b} {\sqrt{6}}},{-\tfrac {b} {\sqrt{2}} } \bigr ),\bigl ({\tfrac {a} {\sqrt{6}}},-{\tfrac {a} {\sqrt{2}} } \bigr )  \bigr \}.$$
The above hexagon is convex when $b<2a$, concave when $b>2a$, and degenerates to an equilateral triangle when $b=2a$   (see Figure~\ref{Fig-E[a,b]}).

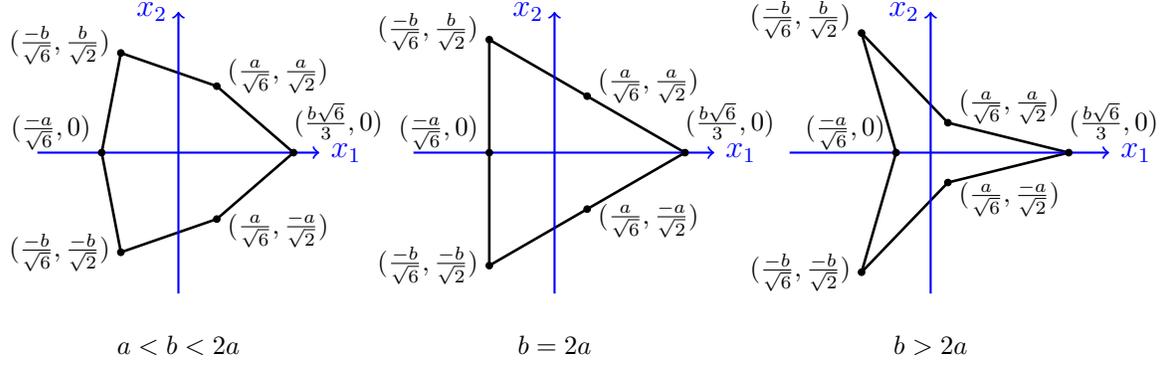
\begin{figure}[t] 
\begin{center}
  \begin{tikzpicture}[scale=1.25]

\path coordinate (q1) at (1.5*0.816,0) coordinate (q2) at (0.408,0.707) coordinate (q3) at  (-1.5*0.408,1.5*0.707) coordinate (q4) at (-0.816,0) coordinate (q5) at (-1.5*0.408,-1.5*0.707) coordinate (q6) at (0.408,-0.707)  ;

\draw[line width=0.35mm] (q1) -- (q2);
\draw[line width=0.35mm] (q2) -- (q3);
\draw[line width=0.35mm] (q3) -- (q4);
\draw[line width=0.35mm] (q4) -- (q5);
\draw[line width=0.35mm] (q5) -- (q6);
\draw[line width=0.35mm] (q6) -- (q1);

\draw[thick,->,blue] (-1.5,0)--(1.5,0) node[right] {$x_1$}; 
    \draw[thick,->,blue] (0,-1.5)--(0,1.5) node[left] {$x_2$}; 

\filldraw [black]
(0,-1.5) circle(0pt) node[below=12pt] {{\footnotesize $a<b<2a$}}

(q1) circle (1pt) node[above=12pt, right=-4pt] {{\footnotesize $({\frac {b 
\sqrt{6}} {3} },0)$}}
      (q2) circle (1pt) node[above=4pt, right] {{\footnotesize $({\frac {a} {\sqrt 6} },{\frac {a} { \sqrt{2}}})$}}
         (q3) circle (1pt) node[above=3pt, left] {{\footnotesize $({\frac {-b} {\sqrt 6} },{\frac {b} { \sqrt{2}}})$}}
          (q4) circle (1pt) node[above=8pt, left] {{\footnotesize $({\frac {-a} {\sqrt 6} },0)$}}
           (q5) circle (1pt) node[ left] {{\footnotesize $({\frac {-b} {\sqrt 6} },{\frac {-b} {\sqrt{2}} })$}}
             (q6) circle (1pt) node[ below=4pt, right] {{\footnotesize $({\frac {a} {\sqrt 6} },{\frac {-a} {\sqrt{2}} })$}};

      \begin{scope}[xshift=4cm]

\path coordinate (p1) at (0.85*2*0.816,0.85*0) coordinate (p2) at (0.85*0.408,0.85*0.707) coordinate (p3) at  (-0.85*2*0.408,0.85*2*0.707) coordinate (p4) at (-0.85*0.816,0.85*0) coordinate (p5) at (-0.85*2*0.408,-0.85*2*0.707) coordinate (p6) at (0.85*0.408,-0.85*0.707)  ;

\draw[line width=0.35mm] (p1) -- (p2);
\draw[line width=0.35mm] (p2) -- (p3);
\draw[line width=0.35mm] (p3) -- (p4);
\draw[line width=0.35mm] (p4) -- (p5);
\draw[line width=0.35mm] (p5) -- (p6);
\draw[line width=0.35mm] (p6) -- (p1);

\draw[thick,->,blue] (-1.5,0)--(1.7,0) node[right] {$x_1$}; 
    \draw[thick,->,blue] (0,-1.5)--(0,1.5) node[left] {$x_2$}; 

\filldraw [black]
(0,-1.5) circle(0pt) node[below=12pt] {{\footnotesize $b=2a$}}
(p1) circle (1pt) node[above=12pt, right=-4pt] {{\footnotesize $({\frac {b 
\sqrt{6}} {3} },0)$}}
      (p2) circle (1pt) node[above = 4pt, right] {{\footnotesize $({\frac {a} {\sqrt 6} },{\frac {a} { \sqrt{2}}})$}}
         (p3) circle (1pt) node[above=3pt, left] {{\footnotesize $({\frac {-b} {\sqrt 6} },{\frac {b} { \sqrt{2}}})$}}
          (p4) circle (1pt) node[above=8pt, left] {{\footnotesize $({\frac {-a} {\sqrt 6} },0)$}}
           (p5) circle (1pt) node[ left] {{\footnotesize $({\frac {-b} {\sqrt 6} },{\frac {-b} {\sqrt{2}} })$}}
             (p6) circle (1pt) node[ below=4pt, right] {{\footnotesize $({\frac {a} {\sqrt 6} },{\frac {-a} {\sqrt{2}} })$}};


\end{scope}

  \begin{scope}[xshift=8cm]

\path coordinate (qq1) at (0.45*4*0.816,0.45*0) coordinate (qq2) at (0.45*0.408,0.45*0.707) coordinate (qq3) at  (-0.45*4*0.408,0.45*4*0.707) coordinate (qq4) at (-0.45*0.816,0.45*0) coordinate (qq5) at (-0.45*4*0.408,-0.45*4*0.707) coordinate (qq6) at (0.45*0.408,-0.45*0.707)  ;

\draw[line width=0.35mm] (qq1) -- (qq2);
\draw[line width=0.35mm] (qq2) -- (qq3);
\draw[line width=0.35mm] (qq3) -- (qq4);
\draw[line width=0.35mm] (qq4) -- (qq5);
\draw[line width=0.35mm] (qq5) -- (qq6);
\draw[line width=0.35mm] (qq6) -- (qq1);

\draw[thick,->,blue] (-1.5,0)--(1.9,0) node[right] {$x_1$}; 
    \draw[thick,->,blue] (0,-1.5)--(0,1.5) node[left] {$x_2$}; 

\filldraw [black]
(0,-1.5) circle(0pt) node[below=12pt] {{\footnotesize {$b>2a$}}}
(qq1) circle (1pt) node[above=12pt, right=-4pt] {{\footnotesize $({\frac {b 
\sqrt{6}} {3} },0)$}}
      (qq2) circle (1pt) node[above=6pt, right] {{\footnotesize $({\frac {a} {\sqrt 6} },{\frac {a} { \sqrt{2}}})$}}
         (qq3) circle (1pt) node[above=5pt, left] {{\footnotesize $({\frac {-b} {\sqrt 6} },{\frac {b} { \sqrt{2}}})$}}
          (qq4) circle (1pt) node[above=8pt, left] {{\footnotesize $({\frac {-a} {\sqrt 6} },0)$}}
           (qq5) circle (1pt) node[ left] {{\footnotesize $({\frac {-b} {\sqrt 6} },{\frac {-b} {\sqrt{2}} })$}}
             (qq6) circle (1pt) node[ below=6pt, right] {{\footnotesize $({\frac {a} {\sqrt 6} },{\frac {-a} {\sqrt{2}} })$}};

\end{scope}
\end{tikzpicture}

\caption{The (possibly degenerate) hexagons that are the closure of $ \mathfrak{P}_E^2[a,b]$ for various values of $a,b >0$.} \label{Fig-E[a,b]}

\end{center}
\end{figure}

By construction, $ \mathfrak{P}_E^2[a,b]$ is $S_3$-invariant, open, bounded  and star-shaped; moreover, its image under the map $\rho:\R^{3}\to\R^{n}_{\ge 0}$  defined in~(\ref{eqn-def-rho}) is the open triangle with vertices $(0, 0), (a, 0)$, and $(0, b)$. Hence, the following is an immediate consequence of Theorem~\ref{thm-lag-prod-as-toric-domain}.
\begin{cor}\label{cor:ellipsoid}
The Lagrangian product  $\mathfrak{S}^2 \times_L \mathfrak{P}_E^2[a,b]$ is symplectically equivalent to $E(3a,3b)$.
\end{cor} 
We can strengthen Corollary \ref{cor:ellipsoid} to the following result.

\begin{cor}
    For all $a,b>0$, $\mathfrak S^2 \times_L \mathfrak{P}_E^2[a,b]$ is symplectomorphic to $E(3a,3b)$.
\end{cor}
\begin{proof}
    The argument is identical to that of \cite[Corollary 4.2 and Lemma 4.3]{LMS} with the only difference that we use \cite[Corollary 1.6 (i)]{McDuff} instead of \cite{Mc}.
\end{proof}

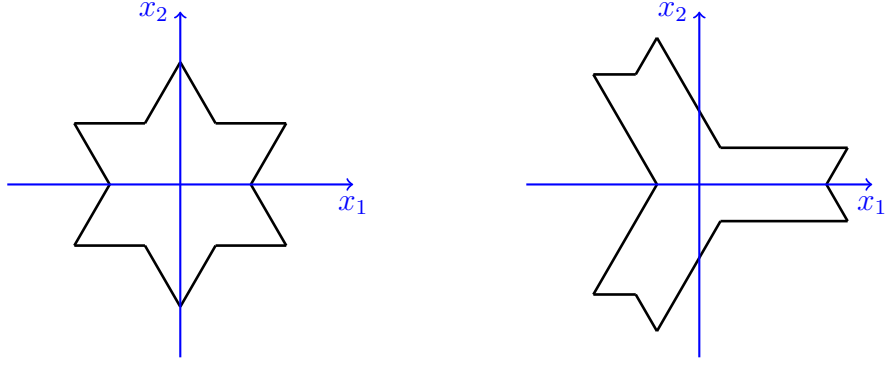
\begin{figure}[h]

\begin{center}
  \begin{tikzpicture}[scale=1.15]
\path coordinate (q1) at (0.816,0)
coordinate (q2) at (1.224,0.707)
coordinate (q3) at (0.408,0.707)
coordinate (q4) at (0,1.4142)
coordinate (q5) at (-0.4082,0.707)
coordinate (q6) at (-1.2247,0.707)
coordinate (q7) at (-0.816,0)
coordinate (q8) at (-1.2247,-0.707)
coordinate (q9) at (-0.4082,-0.707)
coordinate (q10) at (0,-1.4142)
coordinate (q11) at (0.408,-0.707)
coordinate (q12) at (1.22474,-0.707);


\draw[line width=0.35mm] (q1) -- (q2);
\draw[line width=0.35mm] (q2) -- (q3);
\draw[line width=0.35mm] (q3) -- (q4);
\draw[line width=0.35mm] (q4) -- (q5);
\draw[line width=0.35mm] (q5) -- (q6);
\draw[line width=0.35mm] (q6) -- (q7);
\draw[line width=0.35mm] (q7) -- (q8);
\draw[line width=0.35mm] (q8) -- (q9);
\draw[line width=0.35mm] (q9) -- (q10);
\draw[line width=0.35mm] (q10) -- (q11);
\draw[line width=0.35mm] (q11) -- (q12);
\draw[line width=0.35mm] (q12) -- (q1);

\draw[thick,->,blue] (-2,0)--(2,0) node[below] {$x_1$}; 
    \draw[thick,->,blue] (0,-2)--(0,2) node[left] {$x_2$}; 

\filldraw [black]
;
      \begin{scope}[xshift=6cm]
0.4*
\path coordinate (qq1) at (0.6*2.4494,0.6*0)
coordinate (qq2) at (0.6*2.8577,0.6*0.707)
coordinate (qq3) at (0.6*0.408,0.6*0.707)
coordinate (qq4) at (-0.6*0.816,0.6*2.828)
coordinate (qq5) at (-0.6*1.2247,0.6*2.1213)
coordinate (qq6) at (-0.6*2.04124,0.6*2.1213)
coordinate (qq7) at (-0.6*0.816,0.6*0)
coordinate (qq8) at (-0.6*2.0412,-0.6*2.1213)
coordinate (qq9) at (-0.6*1.2247,-0.6*2.1213)
coordinate (qq10) at (-0.6*0.816,-0.6*2.8284)
coordinate (qq11) at (0.6*0.408,-0.6*0.707)
coordinate (qq12) at (0.6*2.8577,-0.6*0.707);

\draw[line width=0.35mm] (qq1) -- (qq2);
\draw[line width=0.35mm] (qq2) -- (qq3);
\draw[line width=0.35mm] (qq3) -- (qq4);
\draw[line width=0.35mm] (qq4) -- (qq5);
\draw[line width=0.35mm] (qq5) -- (qq6);
\draw[line width=0.35mm] (qq6) -- (qq7);
\draw[line width=0.35mm] (qq7) -- (qq8);
\draw[line width=0.35mm] (qq8) -- (qq9);
\draw[line width=0.35mm] (qq9) -- (qq10);
\draw[line width=0.35mm] (qq10) -- (qq11);
\draw[line width=0.35mm] (qq11) -- (qq12);
\draw[line width=0.35mm] (qq12) -- (qq1);

\draw[thick,->,blue] (-2,0)--(2,0) node[below] {$x_1$}; 
    \draw[thick,->,blue] (0,-2)--(0,2) node[left] {$x_2$}; 

\filldraw [black]
;
\end{scope}
\end{tikzpicture}

\caption{The dodecagons that are the closure of $\widetilde{\mathfrak{P}}_P^2[1,1]$ (left)  and $\widetilde{\mathfrak{P}}_P^2[1,3]$  (right).}\label{Fig-polydisk}
\end{center}
\end{figure}

The case of the polydisk $P(a,b)$ is very similar. Define  the set $$ \mathfrak{P}_P^2[a,b] :=  \bigsqcup_{\sigma \in S_3}  \sigma \cdot Q[a,b] \subset W^2,$$ where $Q[a,b]$ is the interior of the parallelogram given by $$Q[a,b] := {\rm Conv} \bigl \{ (0,0,0),  ({\tfrac {-2a} 3},{\tfrac a 3},{\tfrac a 3}),({\tfrac {-b} 3},{\tfrac {-b} 3},{\tfrac {2b} 3}), ( {\tfrac {-2a-b} {3} },{\tfrac {s-b} {3} }, {\tfrac {a+2b} {3} }     )     \bigr \}.$$
Denote $\widetilde{\mathfrak{P}}_P^2[a,b] = (A^{T})^{-1} \mathfrak{P}^2_P[a,b]$,
where $A \in {\rm O}(3)$ as above. The closure of $\widetilde{\mathfrak{P}}_P^2[a,b]$ is the dodecagon composed of the six parallelograms obtained from $Q[a,b]$ by the action of $S_3$ (see Figure~\ref{Fig-polydisk}) with a counterclockwise sequence of vertices given by
\begin{equation*}
\begin{split}
\Bigl \{ & \bigl ({\tfrac {b\sqrt{6}} {3} },0 \bigr ), ( {\tfrac {b\sqrt{6}} {3} } + {\tfrac {a} {\sqrt{6}} },{\tfrac {a} {\sqrt{2}}} ) ,  \bigl ({\tfrac {a} {\sqrt{6}}},{\tfrac {a} {\sqrt{2}} } \bigr ) , ( {\tfrac {a-b} {\sqrt{6}}} ,   {\tfrac {a+b} {\sqrt{2}}}) , \bigl ({\tfrac {-b} {\sqrt{6}}},{\tfrac {b} {\sqrt{2}} } \bigr ) , \bigl ( -{\tfrac {b} {\sqrt{6}}} - {\tfrac {a\sqrt{6}} {3} } , {\tfrac {b} {\sqrt{2}}} \bigr ),  \bigl (-{\tfrac {a\sqrt{6}} {3} },0 \bigr )\\
& \bigl ( -{\tfrac {b} {\sqrt{6}}} - {\tfrac {a\sqrt{6}} {3} } , {-\tfrac {b} {\sqrt{2}}} \bigr ),  \bigl (-{\tfrac {b} {\sqrt{6}}},{-\tfrac {b} {\sqrt{2}} } \bigr ),  \bigl (  {\tfrac {a-b} {\sqrt{6}} }, {\tfrac {-a-b} {\sqrt{2}} }  \bigr ),  \bigl ({\tfrac {a} {\sqrt{6}}},-{\tfrac {a} {\sqrt{2}} } \bigr ),  \bigl (  {\tfrac {b\sqrt{6}} {3} } + {\tfrac {a} {\sqrt{6}} } ,  -{\tfrac {a} {\sqrt{2}} }\bigr ) \Bigr \}. 
\end{split}
\end{equation*}

The set $ \mathfrak{P}_P^2[a,b]$ is $S_3$-invariant, open, bounded, star-shaped and its image under  $\rho$  
 is the open rectangle with vertices $(0, 0), (a, 0), (b,0)$, and $(a, b)$. Hence, the following result is an immediate consequence of Theorem~\ref{thm-lag-prod-as-toric-domain}.
\begin{cor}
The Lagrangian product  $\mathfrak{S}^2 \times_L \mathfrak{P}_P^2[a,b]$ is symplectically equivalent to the polydisk $P(3a,3b)$ for every $a,b>0$.
\end{cor}

\section{The periodic Toda lattice} \label{section-Toda}

\subsection{Integrability via the Lax pair formalism}\label{sec:lax}
The periodic Toda lattice~\cite{To1,To2} is a well-known Hamiltonian system on $\R^{2n}$ describing the movement of $n$ labeled particles on a circle with 
exponential interactions  between nearest neighbors. This is one of the earliest examples of a completely integrable non-linear Hamiltonian system~\cite{Fl1,Hen, Man}, and has been studied
extensively. Before describing it explicitly, we set all indices to be periodic modulo $n$ unless otherwise stated. (We observe that this is a different convention to that of Sections \ref{sec-int} and \ref{sec-examples}.) The Hamiltonian function determining the equations of motion of the periodic Toda lattice is
\begin{equation*}
H(q_1,\dots,q_n,p_1,\dots,p_n)=\frac{1}{2}\sum_{i=1}^n p_i^2+\alpha^2 \sum\limits_{i=1}^n e^{q_i-q_{i+1}},
\end{equation*}
where $\alpha>0$ is a fixed parameter that can be considered as the strength of the potential.
Its equations of motion are given by:
\begin{equation*}
\left\{\begin{aligned}
    \dot{q}_i&=p_i\\
    \dot{p}_i&=\alpha^2\left(e^{q_{i-1}-q_i}-e^{q_i-q_{i+1}}\right),
\end{aligned}\right.
\end{equation*}
where we set $q_{n+1}=q_1$ and $q_0 = q_n$. It is well-known that the periodic Toda lattice is completely integrable (see, for instance, \cite{Fl1,Fl2,FM,vM} among others). One way to construct the desired integrals of motion is to start with a `change of variables' due to Flaschka~\cite{Fl1,Fl2} that describes the motion relative to the center of mass. It is given by the map
\begin{equation}\label{eq:varphialpha}
\begin{split}
    \widetilde{\varphi}_{\alpha} : \R^{2n} &\to \R^n \times \R^n_{>0} \\
    (q_1,\ldots, q_n,p_1,\ldots, p_n) &\mapsto (p_1,\dots,p_n,\alpha e^{\frac{q_1-q_2}{2}},\dots,\alpha e^{\frac{q_n-q_1}{2}}).
\end{split}\end{equation}
We remark that $\widetilde{\varphi}_{\alpha}$ is {\em not} a diffeomorphism: for instance, the image of $\widetilde{\varphi}_{\alpha}$ is contained in the subset 
$$\left\{(b_1,\ldots,b_n,a_1,\ldots, a_n) \in \R^n \times \R^n_{>0} \mid \prod\limits_{i=1}^n a_i = \alpha^n \right\}. $$
If we set $\pi_0 := \omega^{-1}_0$, i.e., $\pi_0 = \sum\limits_{i=1}^n \frac{\partial }{\partial p_i} \wedge \frac{\partial }{\partial q_i}$, then $\pi_0$ can be pushed forward to a bivector field $\pi$ on $\R^n \times \R^n_{>0}$. Explicitly, 
\begin{equation*}
    \pi=\sum_{i=1}^n \frac{a_i}{2}\frac{\partial}{\partial b_i}\wedge \frac{\partial}{\partial a_i}-\sum_{i=1}^n \frac{a_{i-1}}{2}\frac{\partial}{\partial b_i}\wedge \frac{\partial}{\partial a_{i-1}}.
\end{equation*}
Since $\pi_0$ is Poisson, by construction $\pi$ is also Poisson and the map $\tilde{\varphi}_{\alpha}$ is a Poisson morphism. For our purposes, it is useful to prove the following well-known properties of the Poisson manifold $(\R^n \times \R^n_{>0},\pi)$.

\begin{lem}\label{lemma:poisson}
The Poisson structure $\pi$ on $\R^n \times \R^n_{>0}$ has constant rank equal to $2n-2$ and the symplectic leaves are given by the level sets of the map
\begin{equation}
\label{eq:Casimir}
    \begin{split}
        \R^n \times \R^n_{>0} &\to \R^2 \\
        (b_1,\ldots, b_n, a_1,\ldots, a_n) &\mapsto \left(\sum\limits_{i=1}^n b_i \, , \, \prod\limits_{i=1}^n a_i\right).
    \end{split}
\end{equation}
\end{lem}

\begin{proof}
    Let $\{\cdot,\cdot\}$ denote the Poisson bracket on $C^{\infty} (\R^n \times \R^n_{>0})$ induced by $\pi$. For any $i,j = 1,\ldots, n$, we have that $\{b_i,b_j\} = \{a_i,a_j\} = 0$; moreover, $\{b_i,a_j\} = 0$ unless $j=i-1$ or $j=i$. Since the matrix $(\{b_i,a_j\})$ has rank $n-1$, it follows that $\pi$ has constant rank equal to $2n-2$. By a direct computation, the components of the map \eqref{eq:Casimir} are Casimirs of $\pi$, i.e., they have vanishing Hamiltonian vector field. Since the above map is a submersion with connected level sets, and since the rank of $\pi$ is equal to $2n-2$, the result  follows.
\end{proof}

By Lemma \ref{lemma:poisson}, the submanifold
\[\mathcal{M}_{\alpha} :=\left\{(b_1,\dots,b_n,a_1,\dots,a_n)\in\R^n\times\R^n_{> 0}\mid \sum_{i=1}^n b_i=0\text{ and }\prod_{i=1}^n a_i=\alpha^n\right\} \subset \R^n \times \R^n_{>0}\]
is a leaf in the symplectic foliation of $(\R^n \times \R^n_{>0},\pi)$ for each $\alpha >0$. We denote the symplectic form on $\mathcal{M}_{\alpha}$ by $\omega_{\alpha}$. By a direct calculation, $\omega_{\alpha}$ is the restriction to $\mathcal{M}_{\alpha}$ of the 2-form
\begin{equation*}
 2\sum_{i=2}^n \sum_{j=1}^{i-1} db_i\wedge \frac{da_j}{a_j} \in \Omega^2(\R^n \times \R^n_{>0}).
\end{equation*}

Next, there exists a function $\overline{H} \in C^{\infty} (\R^n \times \R^n_{>0})$ such that $\tilde{\varphi}^*_{\alpha} \overline{H} = H$, namely, 
\[ \overline{H}(b_1,\dots,b_n,a_1,\dots,a_n)=\frac{1}{2}\sum_{i=1}^n b_i^2+\sum_{i=1}^n a_i^2.\]
The corresponding equations of motion with respect to $\pi$ are
\begin{equation}\label{eq:todafla}
    \left\{\begin{aligned}
        \dot{b}_i&=a_{i-1}^2-a_{i}^2,\\
        \dot{a}_i&=\frac{a_i}{2}(b_{i}-b_{i+1}).
    \end{aligned}\right.
\end{equation}
The key observation from \cite{Fl1} is that \eqref{eq:todafla} admits a description in terms of a {\bf Lax pair}, i.e., a pair of $n \times n$ matrices $(L,B)$ such that \eqref{eq:todafla} is equivalent to 
\begin{equation}\label{eq:lax}
    \dot{L}=[L,B].
\end{equation}
Explicitly, 
\begin{equation}\label{eq:lb}
L=\begin{bmatrix}
b_1&a_1&0&\dots&0&a_n\\
a_1&b_2&a_2&\dots&0&0\\
0&a_2&b_3&\ddots&0&0\\
\vdots&\vdots&\ddots&\ddots&\ddots&\vdots\\
0&0&0&\ddots&\ddots&a_{n-1}\\
a_n&0&0&\dots&a_{n-1}&b_n
\end{bmatrix},\,
B=\frac{1}{2}\begin{bmatrix}
0&a_1&0&\dots&0&-a_n\\
-a_1&0&a_2&\dots&0&0\\
0&-a_2&0&\ddots&0&0\\
\vdots&\vdots&\ddots&\ddots&\ddots&\vdots\\
0&0&0&\ddots&\ddots&a_{n-1}\\
a_n&0&0&\dots&-a_{n-1}&0
\end{bmatrix}.
\end{equation}
By \eqref{eq:lax}, the eigenvalues $\lambda_1^+\le \dots\le \lambda_n^+$ of $L$ are conserved for the periodic Toda lattice, and so is $\mathrm{Tr}(L) = \sum\limits_{i=1}^n b_i$. If we set $W^{n-1}_{+}=\{(y_1,\dots,y_n)\in W^{n-1}\mid y_1\le \dots \le y_n\}$, the map $F_\alpha:\mathcal{M}_\alpha\to W^{n-1}_{+}$
\begin{equation}
    \label{eq:falpha}
F_\alpha(b_1,\dots,b_n,a_1,\dots,a_n)=(\lambda_1^+,\dots,\lambda_n^+) 
\end{equation}
is a completely integrable Hamiltonian system\footnote{Strictly speaking, the function $F_{\alpha}$ is only smooth on the open dense subset of $\mathcal{M}_{\alpha}$ consisting of points $(b_1,\dots,b_n,a_1,\dots,a_n)$ such that $(\lambda_1^+,\dots,\lambda_n^+)$ are distinct. The map $\tilde{F}_{\alpha} \colon \mathcal{M}_{\alpha} \to W^{n-1}$ with components given by the coefficients of the characteristic polynomial of the matrix $L$ as in \eqref{eq:lb} is smooth and, hence, a completely integrable Hamiltonian system on $(\mathcal{M}_{\alpha},\omega_\alpha)$. Since there is a homeomorphism $i_{\alpha} : \tilde{F}_{\alpha}(\mathcal{M}_{\alpha}) \to W^{n-1}_+$ such that $F_{\alpha} = i_{\alpha} \circ \tilde{F}_{\alpha}$, we abuse terminology slightly and refer to $F_{\alpha}$ as a completely integrable system. We trust that this does not cause confusion.} on $(\mathcal{M}_{\alpha},\omega_\alpha)$ for each $\alpha >0$ (see, e.g., \cite[Section 3.4]{MZ}).

\begin{rem}\label{rmk:beyond_An} {\rm
    As first observed by Bogoyavlensky in~\cite{Bog}, starting from an irreducible root lattice (or, equivalently, a simple Lie algebra), it is possible to construct a Hamiltonian system that fits the Lax formalism. The periodic Toda lattice is precisely the system that arises in the case of the $A_{n-1}$ lattice. Moreover, any such 
 a system can be shown to be integrable using arguments analogous to the ones above (see, for instance, \cite[Section 9.1]{avmm}). }
\end{rem}

\subsection{Global action variables for the periodic Toda lattice}

In this section, we construct global action variables for the periodic Toda lattice, following \cite{FM,HeKa,vM}. 
More precisely, our aim is to construct a symplectomorphism $\psi_{\alpha}: (\mathcal{M}_{\alpha},\omega_\alpha) \to \R^{2n-2}$ and a diffeomorphism $\mathbf{I}_\alpha:W^{n-1}_{+}\to\R^{n-1}_{\ge 0}$ such that $\mu \circ \psi_{\alpha}=\mathbf{I}_\alpha\circ F_\alpha$,
where $\mu:\R^{2n-2}= \C^{n-1}\to \R^{n-1}_{\ge 0}$ is the moment map of the standard toric action on $\R^{2n-2} \simeq \C^{n-1}$, i.e.,
\begin{equation}
    \label{eq:mu}
    \mu(z_1,\dots,z_{n-1})=(\pi|z_1|^2,\dots,\pi |z_{n-1}|^2).
\end{equation}
The components of the map $\mathbf{I}_{\alpha}$ are the desired global action variables.

In order to proceed, we recall some facts about a discrete version of Floquet theory (for more details and proofs, see \cite{vM}). To this end, we define $a_k$ and $b_k$ for all $k\in\mathbb{Z}$ by letting $a_{k+n}=a_k$ and $b_{k+n}=b_k$, and consider the difference equation 
\begin{equation}
    \label{eq:diff}a_{k-1}y(k-1,\lambda)+b_k y(k,\lambda)+a_k y(k+1,\lambda)=\lambda y(k,\lambda), \quad k=1,2,\dots,
\end{equation}
for a given $\lambda \in \R$.
An eigenvector of the matrix $L$ of \eqref{eq:lb} with eigenvalue $\lambda$ is precisely a periodic solution $y(k,\lambda)$ of \eqref{eq:diff}, i.e., $y(k+n,\lambda)=y(k,\lambda)$. We define two fundamental solutions $y_1(k,\lambda)$ and $y_2(k,\lambda)$ of \eqref{eq:diff} by letting
\begin{equation}\label{eq:init}
    y_1(0,\lambda)=1, \quad y_1(1,\lambda)=0,\quad y_2(0,\lambda)=0,\quad y_2(1,\lambda)=1.
\end{equation}
Any solution to \eqref{eq:diff} for a fixed $\lambda$ is a linear combination of $y_1(k,\lambda)$ and $y_2(k,\lambda)$. By a simple calculation,
\[y_1(k,\lambda)y_2(k+1,\lambda)- y_1(k+1,\lambda)y_2(k,\lambda)=\frac{a_{k-1}}{a_k}(y_1(k-1,\lambda)y_2(k,\lambda)-y_1(k,\lambda)y_2(k-1,\lambda)),\]
for any $k$. Hence, since $a_{n}=a_0$,
\begin{equation}
    \label{eq:wr}y_1(n,\lambda)y_2(n+1,\lambda)- y_1(n+1,\lambda)y_2(n,\lambda)=1.
\end{equation}
The solution $y_1(n+1,\lambda)$ is a polynomial of degree $n-1$ in $\lambda$. As explained in \cite[Section 2]{vM}, the roots of $y_1(n+1,\lambda)$ are precisely the eigenvalues of the symmetric matrix obtained by removing the first row and column from the matrix $L$. In particular, they are real; we denote them by  $\mu_1\le \dots\le \mu_{n-1}$. By \eqref{eq:wr}, $y_2(n+1,\mu_i)\neq 0$ for any $i=1,\dots,n-1$. Hence, the quantity  
$$f_i: =\log|y_2(n+1,\mu_i)|$$ 
is well-defined for any $i=1,\dots,n-1$. The following result is essentially proved in \cite{vM}, and claimed in \cite{FM}, for $\alpha=1$. For the sake of completeness, we give a proof.

\begin{thm}\label{thm:phisymp}
The image of the map $\phi_{\alpha}:\mathcal{M}_\alpha\to \R^{2n-2} $ that sends $(b_1,\dots,b_n,a_1,\dots,a_n)$ to $(f_1,\dots,f_{n-1},\mu_1,\dots,\mu_{n-1})$ consists of all values such that $\mu_1<\mu_2<\dots<\mu_{n-1}$. Moreover, $\phi_{\alpha}$ is a symplectomorphism from $(\mathcal{M}_\alpha,\omega_{\alpha})$ onto $(\phi_{\alpha}(\mathcal{M}_\alpha), \omega_0=\sum_{i=1}^{n-1} d\mu_i\wedge d f_i)$.
In particular, the 1-form $\phi_{\alpha}^*\left(\sum_{i=1}^{n-1} f_i \,d\mu_i\right)$ is a primitive of $\omega_{\alpha}$.
\end{thm}

\begin{proof}
We fix $(b_1,\dots,b_n,a_1,\dots,a_n)\in\mathcal{M}_\alpha$ and let 
$\mu_1\le\dots\le \mu_{n-1}$ be as above. First, we show that $\mu_1<\mu_2<\dots<\mu_{n-1}$.
By \eqref{eq:diff}, 
for any $\lambda,\widetilde{\lambda} \in \R$, 
\begin{equation*}
    (\lambda-\widetilde{\lambda})\sum_{k=2}^n y_1(k,\lambda)y_1(k,\widetilde{\lambda})=a_n\left(y_1(n+1,\lambda)y_1(n,\widetilde{\lambda})-y_1(n+1,\widetilde{\lambda})y_1(n,\lambda)\right).
\end{equation*}
Hence, 
\begin{equation}\label{eq:vect}
    \sum_{k=2}^n y_1(k,\mu_i)y_1(k,\mu_j)=
    \begin{cases}
    a_n y_1'(n+1,\mu_i)y_1(n,\mu_i) &\text{if } i=j,\\
    0 &\text{if }i\neq j, 
    \end{cases}
\end{equation}
where $y_1'(n+1,\mu_i)$ is the derivative of $y_1(n+1,\lambda)$ with respect to $\lambda$ evaluated at $\lambda=\mu_i$.
For each $i=1,\dots,n-1$, let $v_i\in\R^{n-1}$ be the vector given by $(y_1(2,\mu_i),\ldots, y_1(n,\mu_i))$. By \eqref{eq:vect}, $v_i$ and $v_j$ are orthogonal whenever $i \neq j$. On the other hand, by \eqref{eq:diff}, the leading coefficient of $y_1(n+1,\lambda)$ is $-a_0/\alpha^n=-a_n/\alpha^n$. Hence,
\begin{equation}\label{eq:y1}
y_1(n+1,\lambda)=-\frac{a_n}{\alpha^n}\prod_{i=1}^{n-1}(\lambda-\mu_i).
\end{equation}
By \eqref{eq:wr}, \eqref{eq:vect}, and \eqref{eq:y1}, 
\begin{equation}\label{eq:norm}\Vert v_i\Vert^2=a_n y_1'(n+1,\mu_i)y_1(n,\mu_i)=-\frac{a_n^2\prod_{j\neq i}(\mu_i-\mu_j)}{\alpha^n\cdot  y_2(n+1,\mu_i)} = \frac{a_n^2 \prod_{j\neq i}|\mu_i - \mu_j|}{\alpha^n e^{f_i}},\end{equation} 
where the last equality follows from the fact that $\Vert v_i \Vert \geq 0$. Moreover, by \eqref{eq:diff}, $y_1(2,\mu_i)=-a_n/a_1$, so that $\Vert v_i\Vert\neq 0$ for any $i = 1,\ldots, n-1$. Hence, by \eqref{eq:norm}, $\mu_i\neq\mu_j$ for all $i\neq j$ and, in particular, $\mu_1<\dots<\mu_{n-1}$. As an immediate consequence, $\phi_{\alpha}$ is smooth.

Next we show that $\phi_{\alpha}$ is injective and compute its image. We set \[(\mathbf{f},\boldsymbol{\mu})=(f_1,\dots,f_{n-1},\mu_1,\dots,\mu_{n-1}):=\phi_{\alpha}(b_1,\dots,b_n,a_1,\dots,a_n).\]
Since $\{v_i/\Vert v_i\Vert\}_{i=1,\dots,n-1}$ is an orthonormal basis of $\R^{n-1}$, 

\begin{equation}\label{eq:on}\sum_{i=1}^{n-1} \frac{y_1(k,\mu_i)y_1(j,\mu_i)}{\Vert v_i\Vert^2}=\delta_{kj},\end{equation} for $2\le k,j\le n$, where $\delta_{kj}$ is the Kronecker delta. For $i=1,\dots,n-1$, we introduce weights 
\begin{equation}
\label{eq:weights}
 w_i:= a_n^2/\Vert v_i\Vert^2 = \frac{\alpha^n e^{f_i}}{\prod_{j\neq i}|\mu_i - \mu_j|}. 
\end{equation}
We define the following bilinear, symmetric product on the polynomial ring $\R[\lambda]$: Given $\mathfrak{p}_1(\lambda),\mathfrak{p}_2(\lambda) \in \R[\lambda]$,
\begin{equation}\label{eq:ip}\left\langle\mathfrak{p}_1(\lambda),\mathfrak{p}_2(\lambda)\right\rangle:=\sum_{i=1}^{n-1}\mathfrak{p}_1(\mu_i)\mathfrak{p}_2(\mu_i)w_i.\end{equation}
By \eqref{eq:weights}, and since the weights are positive, this product is positive semi-definite and depends only on $(\mathbf{f},\boldsymbol{\mu})$. Moreover, $\langle \mathfrak{p}(\lambda), \mathfrak{p}(\lambda) \rangle = 0$ if and only if $\mathfrak{p}(\lambda)$ lies in the ideal generated by $y_1(n+1,\lambda)$. Hence, since the degree of $y_1(n+1,\lambda)$ is $n-1$, the restriction of $\langle \cdot, \cdot \rangle$ to the vector subspace $\mathrm{P}^{n-2}$ of polynomials of degree at most $n-2$ is, in fact, an inner product, which we also denote by $\langle \cdot, \cdot \rangle$. If we set $P_k(\lambda):=y_1(k,\lambda)/a_n$ for $k=2,\dots,n$, then $\{P_2(\lambda),\ldots, P_n(\lambda)\}$ is a basis of $\mathrm{P}^{n-2}$ that has the following properties:
\begin{enumerate}[label = (\alph*),ref=(\alph*), leftmargin=*]
    \item \label{item:-1} it is orthonormal by \eqref{eq:on} and \eqref{eq:ip};
    \item \label{item:-2} the degree of $P_k(\lambda)$ is $k-2$ for each $k=2,\ldots, n$; and
    \item \label{item:-3} the sign of the leading term of $P_k(\lambda)$ is $-1$ for each $k=2,\ldots, n$.
\end{enumerate}
Since there is precisely one such basis of $\mathrm{P}^{n-2}$, the polynomials $\{P_k\}_{k=2,\dots,n}$ are uniquely determined by \eqref{eq:ip}, so that they depend only on $(\mathbf{f},\boldsymbol{\mu})$. By \eqref{eq:diff}, property \ref{item:-1} and the definition of $\mathcal{M}_{\alpha}$,
\begin{equation}\label{eq:ab}
    \begin{aligned}
    b_k&=\left\langle  b_k P_k(\lambda),P_k(\lambda)\right\rangle=\langle \lambda P_k(\lambda),P_k(\lambda)\rangle,&k=2,\dots,n,&\\
     b_1&=-(b_2+\dots+b_n),\\
    a_k&=\left\langle  a_k P_{k+1}(\lambda),P_{k+1}(\lambda)\right\rangle=\langle \lambda P_k(\lambda),P_{k+1}(\lambda)\rangle,&k=2,\dots,n-1,&\\
    a_1&=-\frac{1}{P_2(\lambda)},\\
    a_n&=\frac{\alpha^n}{a_1\dots a_{n-1}}.
    \end{aligned}
\end{equation}
Hence, by \eqref{eq:ab}, $\phi_{\alpha}$ is injective. In fact, it is onto the subset of $\R^{2n-2}$ that consists of points $(\mathbf{f},\boldsymbol{\mu})$ such that $\mu_1<\dots<\mu_{n-1}$: To see this, first we use \eqref{eq:weights} to define an inner product on $\mathrm{P}^{n-2}$ as in \eqref{eq:ip}. Second, we consider the unique basis of $\mathrm{P}^{n-2}$ that satisfies properties \ref{item:-1} -- \ref{item:-3} above and use \eqref{eq:ab} to {\em define} $(b_1,\ldots, b_n,a_1,\ldots,a_n)$. This is the desired pre-image of $(\mathbf{f},\boldsymbol{\mu})$ under $\phi_{\alpha}$.

If we show that $\phi_{\alpha}^*\omega_0=\omega_{\alpha}$, then $\phi_{\alpha}$ is an injective immersion between manifolds of equal dimension and, therefore, a symplectomorphism onto its image, as desired. This is proved in \cite[Theorem 7.1]{vM} (where, strictly speaking, only the case $\alpha=1$ is considered, but the proof works {\em mutatis mutandis}).
\end{proof}

In order to construct the map $\mathbf{I}_\alpha$, we recall some more facts from discrete Floquet theory (see \cite{vM}). We set 
\begin{equation}
   \label{eq:delta_def} \Delta(\lambda):=y_1(n,\lambda)+y_2(n+1,\lambda).
\end{equation}
By a simple calculation, the difference equation \eqref{eq:diff} admits a non-trivial solution of period $n$ if and only if $\Delta(\lambda)=2$. Since $\Delta(\lambda)$ and $\det(\lambda I-L)$ are both polynomials of degree $n$, if we denote the eigenvalues of $L$ by $\lambda_1^+ \leq \ldots \leq \lambda_n^+$, then 
\begin{equation}\label{eq:delta}\Delta(\lambda)=\frac{1}{a_1\dots a_n}\det(\lambda I-L)+2=\frac{1}{\alpha^n}\prod_{i=1}^n(\lambda-\lambda_i^+) +2.\end{equation}
Similarly, \eqref{eq:diff} admits a non-trivial  solution satisfying $y(k+n,\lambda)=-y(k,\lambda)$ if and only if $\Delta(\lambda)=-2$. As above, we denote the solutions to $\Delta(\lambda)=- 2$ by $\lambda_1^{-}\le\dots \leq \lambda_n^{-}$. (These are real because they are precisely the eigenvalues of the symmetric matrix obtained by switching the sign of $a_n$ in $L$.) Since non-trivial solutions satisfying $y(k+n,\lambda)=\pm y(k,\lambda)$ have period $2n$, the zeros of $\Delta(\lambda)^2-4$ are precisely the eigenvalues of
\begin{equation}\label{eq:Q}Q:=\begin{bmatrix}b_1&a_1&0&\cdots&0&0&0&0&\cdots&a_{n}\\a_1&b_2&a_2&\cdots&0&0&0&0&\cdots&0\\
0&a_2&b_3&\ddots&0&0&0&0&\cdots&0\\
\vdots&\ddots&\ddots&\ddots&\vdots&\vdots&\ddots&\ddots&\ddots&\vdots\\0&\cdots&\cdots&a_{n-1}&b_{n}&a_{n}&\cdots&0&0&0\\
0&0&0&\cdots&a_{n}&b_1&a_1&0&\cdots&0\\0&0&0&\cdots&0&a_1&b_2&a_2&\cdots&\vdots\\
0&0&0&\cdots&0&0&a_2&b_3&\ddots&0\\\vdots&\ddots&\ddots&\ddots&\vdots&\vdots&\ddots&\ddots&\ddots&\vdots\\a_{n}&\cdots&0&0&0&0&\cdots&\cdots&a_{n-1}&b_{n}
\end{bmatrix}\end{equation} 
We denote the eigenvalues of $Q$ by $\lambda_1\le \lambda_2\le \dots\le \lambda_{2n}$.
As proved in \cite{vM}, $\lambda_{2i-1}<\lambda_{2i}$ for all $i=1,\dots,n-1$, $\lambda_{2n}=\lambda_n^+$, and $\lambda_1=\lambda_{1}^{+}$ or $\lambda_1=\lambda_1^-$, where the sign depends on the parity of $n$. Recall that $\mu_1 < \ldots <\mu_{n-1}$ are the roots of the polynomial $y_1(n+1,\lambda)$. By \eqref{eq:wr} and \eqref{eq:delta_def}, 
\begin{equation}\label{eq:deltamu}
    \Delta(\mu_i)=y_2(n+1,\mu_i)+\frac{1}{y_2(n+1,\mu_i)} \quad \text{for all } i=1,\dots,n-1.
\end{equation}
In particular, 
\begin{equation}
    \label{eq:deltageq2}
    |\Delta(\mu_i)|\ge 2.   
\end{equation}
Hence, by \eqref{eq:y1}, each interval $[\lambda_{2j},\lambda_{2j+1}]$ contains a unique $\mu_i$, so that $\mu_i\in[\lambda_{2i},\lambda_{2i+1}]$ for each $i=1,\dots,n-1$. Moreover, $\sum_i \lambda_i^+=\mathrm{Tr}(L)=0$.

Fix $\lambda_1^+ \leq \dots \leq \lambda_n^+$ satisfying $\sum_i \lambda_i^+=0$. By \eqref{eq:delta}, this choice together with $\alpha$ determine $\Delta(\lambda)$ and hence the spectrum of $Q$. The regular values of the map $F_\alpha$ of \eqref{eq:falpha} are precisely those $\lambda_1^+<\dots<\lambda_n^+$ for which the induced spectrum of $Q$ satisfies $\lambda_{2i}<\lambda_{2i+1}$ for all $i=1,\dots,n-1$. In this case, $F_\alpha^{-1}(\lambda_1^+,\dots,\lambda_n^+)$ is an $(n-1)$-torus in $\mathcal{M}_\alpha$. To see this, we let $\mu_i\in[\lambda_{2i},\lambda_{2i+1}]$ for each $i=1,\ldots, n-1$. By \eqref{eq:deltamu}, there are two possible values of $y_2(n+1,\mu_i)$ (and hence of $f_i$) whenever $\mu_i$ is not equal to either $\lambda_{2i}$ or $\lambda_{2i+1}$ and only one such value when $\mu_i=\lambda_{2i}$ or $\mu_i = \lambda_{2i+1}$. Hence, if we fix a value $\mu_j$ for all $j\neq i$ and let $\mu_i$ vary, we obtain a circle in $F_{\alpha}^{-1}(\lambda_1^+,\dots,\lambda_n^+)$ that we denote by $\gamma_i$. (This circle $\gamma_i$ can be thought of as two copies of $[\lambda_{2i},\lambda_{2i+1}]$ glued at the endpoints.) We orient $\gamma_i$ so that, in the region  $f_i>0$, $\mu_i$ is increasing. By \eqref{eq:deltamu},
\[|f_i|=\Big|\log|y_2(n+1,\mu_i)|\Big|=\Bigg|\log\frac{\big|\Delta(\mu_i)\pm\sqrt{\Delta(\mu_i)^2-4}\big|}{2}\Bigg|=\cosh^{-1}\frac{|\Delta(\mu_i)|}{2}.\]

Let $\mathcal{M}_\alpha^*\subset\mathcal{M}_\alpha$ denote the regular points of the integrable system, i.e., 
$$(b_1,\dots,b_n,a_1,\dots,a_n)\in\mathcal{M}^*_\alpha \Leftrightarrow \lambda_{2i}<\lambda_{2i+1} \quad \text{for all } i=1,\dots,n-1.$$ 
By the standard formula for action variables (see, e.g., \cite[Chapter 10, Section 50, Part C]{Arn}), the $i$-th component of $\mathbf{I}_\alpha$ is given by 
\begin{equation}\label{eq:action}I_i^\alpha(\lambda_1^+,\dots,\lambda_n^+)=\oint_{\gamma_i}\phi^*\left(\sum_{j=1}^{i-1} f_j\, d\mu_j\right)=2\int_{\lambda_{2i}}^{\lambda_{2i+1}} \cosh^{-1}\frac{|\Delta(\mu_i)|}{2}\,d\mu_i.
\end{equation}

{\em A priori}, \eqref{eq:action} only makes sense in a sufficiently small neighborhood of a regular value of $\mathbf{I}_{\alpha}$. However, as proved in \cite[Theorem 1.1]{HeKaaa}, it makes sense on the whole set of regular values $\mathbf{I}_{\alpha}$. Moreover, by \cite[Theorem 1.1]{HeKa}, for each $i=1,\ldots, n-1$, the function of \eqref{eq:action} can be extended to the whole image $F_{\alpha}(\mathcal{M}_{\alpha}) = W^{n-1}_+$ and there exists a symplectomorphism $\psi_{\alpha} : \mathcal{M}_\alpha \to \C^{n-1}$ such that 
\begin{equation}\label{eq:muphi}\mu\circ \psi_\alpha=\mathbf{I}_\alpha\circ F_\alpha,\end{equation}
where $\mu : \C^{n-1} \to \R_{\geq 0}^{n-1}$ is as in \eqref{eq:mu}.


\section{The proof of our main result} \label{section-from-Toda-to-billiards}

\subsection{From the periodic Toda lattice to the billiard in $\mathfrak{S}^{n-1}$}
In this section we introduce the deformation of the periodic Toda lattice 
that we use in the proof of Theorem \ref{thm-lag-prod-as-toric-domain}. Fix $\alpha >0$; moreover, as in Section \ref{section-Toda}, all indices are modulo $n$ unless otherwise stated. Finally, we set $$ (\mathbf{q},\mathbf{p}) = (q_1,\ldots, q_n,p_1,\ldots, p_n) \in \R^{2n}.$$
We recall that $V^{2n-2}:=W^{n-1} \times_L W^{n-1} \subset \R^{2n}$ is the symplectic subspace given by 
$$\left\{(\mathbf{q},\mathbf{p}) \in \R^{2n} \mid \sum_{i=1}^n q_i=0\text{ and } \sum_{i=1}^n p_i=0 \right\}. $$
The following  is an immediate consequence of Lemma \ref{lemma:poisson} and of the definition of $\widetilde{\varphi}_{\alpha}$ in \eqref{eq:varphialpha}. 

\begin{cor}\label{cor:sympleaf}
    For each $\alpha >0$, the restriction of the map $\widetilde{\varphi}_{\alpha} : \R^{2n} \to \R^n \times \R^n_{>0}$ of \eqref{eq:varphialpha} to $V^{2n-2}$ is a symplectomorphism $\varphi_{\alpha} : V^{2n-2} \to (\mathcal{M}_{\alpha},\omega_{\alpha})$.
\end{cor}

Explicitly, the map $\varphi_{\alpha}$ is given by 
\begin{equation}
\label{eq:phialpha}
\varphi_\alpha(\mathbf{q},\mathbf{p})=(p_1,\dots,p_n,\alpha e^{\frac{q_1-q_2}{2}},\dots,\alpha e^{\frac{q_n-q_1}{2}}).
\end{equation} 
We set $\Phi_\alpha: =\psi_\alpha\circ \varphi_\alpha:V^{2n-2}\to \R^{2n-2}$, where $\psi_\alpha$  is the symplectomorphism satisfying $(\ref{eq:muphi})$. By construction, the map $\Phi_{\alpha}$ is a symplectomorphism. 

Next we let $c >0$ and consider the potential $U_c : W^{n-1} \to \R$ given by 
\begin{equation}
\label{eq:uc}
    \mathbf{q} \mapsto c^2e^{-c}\sum_{i=1}^n e^{c(q_i-q_{i+1})}.
\end{equation}
The following result establishes some properties of $U_c$.

\begin{lem}\label{lemma:uc}
    There exists $c_0 >0$ such that the family of subsets $\{U_c^{-1}([0,M])\}_{c \,\geq c_0}$ is monotone increasing for any $M > 0$, i.e.,
    $$ \text{if }c_1 \geq c_2 \geq c_0 \text{ then } U^{-1}_{c_2} ([0,M]) \subseteq U^{-1}_{c_1} ([0,M]).$$
    Moreover, there exists $M_0 >0$ such that, for any $M \geq M_0$.
    $$ \bigcup\limits_{c \,\geq c_0} U_c^{-1}([0,M]) = \mathfrak{S}^{n-1},$$
    \noindent
    where $\mathfrak{S}^{n-1} = \{ \mathbf{q} \in W^{n-1} \mid q_{i} - q_{i+1} < 1 \ {\rm for} \ 1\leq i \leq n\} $ as in \eqref{eq:100}.
\end{lem}

\begin{proof}
    Fix $M >0$. First we show that, given $c >0$, there exists a constant $C>0$ such that 
    $$ U_c^{-1}([0,M]) \subset \{\mathbf{q} \in W^{n-1} \mid q_i - q_{i+1} -1 \leq C \} \text{ for each } i=1,\ldots, n.$$
    If not, then there exist a sequence of points $\{\mathbf{q}^{(k)}\}_{k \in \mathbb{N}} \subset U_c^{-1}([0,M])$ and $i=1,\ldots, n$ such that the sequence $\{q^{(k)}_i - q^{(k)}_{i+1} -1 \}_{k \in \mathbb{N}}$ consists of positive numbers and is unbounded. Hence, the sequence of positive numbers $\{c^2 e^{c(q^{(k)}_i - q^{(k)}_{i+1} - 1)}\}_{k \in \mathbb{N}}$ is unbounded. However, by \eqref{eq:uc}, 
    $$ M \geq U_c(\mathbf{q}^{(k)}) \geq c^2 e^{c(q^{(k)}_i - q^{(k)}_{i+1} - 1)} \text{ for each } k \in \mathbb{N}. $$
    This is a contradiction. 

    Next we show that there exists $c_0 >0$ such that the function that sends $c\geq c_0$ to $U_c(\mathbf{q})$ is decreasing for every $\mathbf{q} \in U^{-1}_c([0,M])$. To this end, we fix $ \mathbf{q}\in U^{-1}_c([0,M])$ and calculate $\partial_c U_c(\mathbf{q})$ for $c \geq c_0$, where $c_0:= 2$:
    \begin{equation*}
        \begin{split}
            \partial_c U_c(\mathbf{q}) &= 2c \sum\limits_{i=1}^n e^{c(q_i - q_{i+1}-1)} + c^2\sum\limits_{i=1}^n (q_i - q_{i+1}-1)e^{c(q_i - q_{i+1}-1)} \\
            & \leq c^2 \sum\limits_{i=1}^n (q_i - q_{i+1})e^{c(q_i - q_{i+1}-1)} \leq \tilde{C} \sum\limits_{i=1}^n (q_i - q_{i+1}) = 0,
        \end{split}
    \end{equation*}
    \noindent
    where $\tilde{C} = c^2e^{cC}$ and we use $c \geq c_0$ in the first inequality. As an immediate consequence, the family of subsets $\{U_c^{-1}([0,M])\}_{c \geq c_0}$ is monotone increasing. Moreover, it also follows that $U^{-1}_c([0,M])$ is contained in $\mathfrak{S}^{n-1}$ for every $c \geq c_0$. To see this, we observe that, given $\mathbf{q} \in W^{n-1}$, 
    \begin{equation}\label{eq:ulim}
    \lim_{c\to + \infty} U_c(\mathbf{q}) =
        \begin{cases}
            0 & \text{if } \mathbf{q} \in \mathfrak{S}^{n-1}, \\
            + \infty & \text{if }\mathbf{q} \in W^{n-1} \smallsetminus \mathfrak{S}^{n-1}.
        \end{cases}
\end{equation}
Since the positive function that sends $c\geq c_0$ to $U_c(\mathbf{q})$ is decreasing, if $\mathbf{q} \in (W^{n-1} \smallsetminus \mathfrak{S}^{n-1}) \cap U^{-1}_{c_1}([0,M])$ for some $c_1 \geq c_0$, then the subset $\{U_c(\mathbf{q})\}_{c \, \geq c_1}$ is bounded above. However, this contradicts \eqref{eq:ulim}. Hence, 
$$ \bigcup\limits_{c \, \geq c_0} U_c^{-1}([0,M]) \subseteq \mathfrak{S}^{n-1}.$$
Conversely, for any $M \geq M_0:= 4n$, if $\mathbf{q} \in \mathfrak{S}^{n-1}$, then
$ U_{c_0}(\mathbf{q}) \leq nc_0^2 = 4n \leq M$,
so that 
$$\mathfrak{S}^{n-1}\subseteq \bigcup\limits_{c \, \geq c_0} U_c^{-1}([0,M]).$$
\end{proof}

\begin{rem}\label{rem:limit}
 Lemma \ref{lemma:uc} implies that the (standard) billiard in the closure of $\mathfrak{S}^{n-1}$ can be seen as a suitable limit of a deformation of the periodic Toda lattice. To see this, we consider the Hamiltonian $H_c : V^{2n-2} \to \R$ given by $H_c(\mathbf{q},\mathbf{p})=\frac{1}{2}|\mathbf{p}|^2+U_c(\mathbf{q})$, where $|\mathbf{p}|^2 = \sum_{i=1}^n p^2_i$ and $U_c(\mathbf{q})$ is as in \eqref{eq:uc}. By Lemma \ref{lemma:uc} and \eqref{eq:ulim}, the limit of the flow of the Hamiltonian vector field $X_{H_c}$ as $c \to + \infty$ is the billiard flow in the closure of $\mathfrak{S}^{n-1}$ (cf. \cite[Chapter IV, Section 5]{KT}).
\end{rem}

\subsection{The proof of Theorem \ref{thm-lag-prod-as-toric-domain}}
Motivated by the above discussion, we define
\begin{equation}\label{eq:psi}
    \begin{split}
        \Psi_c:V^{2n-2} & \to\R^{2n-2} \\
        \Psi_c(\mathbf{q},\mathbf{p}) &=\frac{1}{\sqrt{c}}\Phi_{ce^{-c/2}}(c\, \mathbf{q},\mathbf{p}).
    \end{split}
\end{equation}
By a direct calculation, $\Psi_c$ is a symplectomorphism. We recall that 
$$W^{n-1}_{+}=\{\mathbf{y} :=(y_1,\dots,y_n)\in W^{n-1}\mid y_1\le \dots \le y_n\},$$ 
and define the functions $\mathbf{F}_c:V^{2n-2}\to W^{n-1}_+$ and $\mathbf{J}_c:W^{n-1}_+\to\R_{\geq 0}^{n-1}$ by
\begin{equation}
\begin{aligned}
    \mathbf{F}_c(\mathbf{q},\mathbf{p})&=F_{ce^{-c/2}}\circ \varphi_{ce^{-c/2}}(c \, \mathbf{q},\mathbf{p}),\\    \mathbf{J}_c(\mathbf{y})&=\frac{1}{c}\,\mathbf{I}_c(\mathbf{y}),
    \end{aligned}\label{eq:deffjc}
\end{equation}
where $F_{\alpha} : \mathcal{M}_{\alpha} \to W^+_{n-1}$ is as in \eqref{eq:falpha} and $\mathbf{I}_{\alpha} : W^+_{n-1} \to \R^{n-1}_{\geq 0}$ is the map with components given by \eqref{eq:action}. Let $L_c$ denote the matrix $L$ in \eqref{eq:lb} with $(b_1,\dots,b_n,a_1,\dots,a_n)=\varphi_{ce^{-c/2}}(c\, \mathbf{q},\mathbf{p})$. We set 
\begin{equation}
    \label{eq:lambdaplus}
    (\lambda_1^+(c),\dots, \lambda_n^+(c)):=\mathbf{F}_c(\mathbf{q},\mathbf{p}),
\end{equation}
where $\lambda_1^+(c) \leq \dots \leq \lambda_n^+(c)$ are the eigenvalues of $L_c$. By \eqref{eq:muphi}, \eqref{eq:psi} and \eqref{eq:deffjc}, \begin{equation}\label{eq:alc}\mathbf{J}_c\circ \mathbf{F}_c= \mu\circ\Psi_c.\end{equation}

The following result is the key ingredient in the proof of Theorem \ref{thm-lag-prod-as-toric-domain}.
\begin{prop}\label{prop:lim}
\mbox{}
\begin{enumerate}[label=(\alph*)]
    \item Let $(\mathbf{q},\mathbf{p})\in V^{2n-2}$. If $\mathbf{q}\in\mathfrak S^{n-1}$, then 
    $$\lim_{c\to +\infty}\mathbf{F}_c(\mathbf{q},\mathbf{p})=(p_{\sigma(1)},\dots,p_{\sigma(n)}),$$ 
    where $\sigma \in S_n$ is such that $p_{\sigma(1)}\le \dots\le p_{\sigma(n)}$. Moreover, the convergence is uniform on compact subsets of $\mathfrak S^{n-1}\times_L W^{n-1} \subset V^{2n-2}$. Otherwise, $\lim_{c\to +\infty} |\mathbf{F}_c(\mathbf{q},\mathbf{p})|= +\infty.$
    \item If $\mathbf{y}\in W^{n-1}_+$, then
\begin{equation}
    \label{eq:limJ}
    \lim_{c\to +\infty}\mathbf{J}_c(\mathbf{y})=n(y_2-y_1,\dots,y_n-y_{n-1})=\rho(y_1,\dots,y_n);
\end{equation}
    moreover, the convergence is uniform on compact subsets of $W^{n-1}_+$.
\end{enumerate}

\end{prop}

\begin{proof} \mbox{}
\begin{enumerate}[label=(\alph*)]
\item We observe that, for any $(\mathbf{q},\mathbf{p})\in V^{2n-2}$,
\begin{equation}\label{eq:flim}
|\mathbf{F}_c(\mathbf{q},\mathbf{p})|^2=\sum_{i=1}^n(\lambda_i^+)^2=\mathrm{Tr}(L_c^2)=\sum_{i=1}^n p_i^2+2\sum_{i=1}^nc^2 e^{c(q_i-q_{i+1}-1)}.
\end{equation}
Hence, $|\mathbf{F}_c(\mathbf{q},\mathbf{p})|$ is bounded in $c$ precisely if $\mathbf{q}\in\mathfrak S^{n-1}$, proving the last claim.

If $\mathbf{q}\in\mathfrak S^{n-1}$, then $\lim\limits_{c \to + \infty} L_c$ is a diagonal matrix with entries $(b_1,\dots,b_n)=(p_1,\dots,p_n)$. Hence,
\begin{equation}\label{eq:fclim}\lim_{c\to \infty}\mathbf{F}_c(\mathbf{q},\mathbf{p})=\lim_{c\to\infty} (\lambda_1^+(c),\dots,\lambda_n^+(c))=(p_{\sigma(1)},\dots,p_{\sigma(n)}),\end{equation}
where $\sigma\in S_n$ is such that $p_{\sigma(1)}\le\dots\le p_{\sigma(n)}$. Moreover, if $\mathbf{q}$ belongs to a compact subset of $\mathfrak S^{n-1}$, then there exists $\delta>0$ such that $q_i-q_{i+1}-1<-2\delta$ for each $i=1,\ldots, n$. Hence,  $0<a_i<ce^{-\delta c}$ for every $i$ so that the convergence in \eqref{eq:fclim} is uniform if $\mathbf{p}$ belongs to a fixed compact subset of $W^{n-1}$. 

\item 
Let $ \mathbf{y}\in W^{n-1}_+$ and set $\mathbf{J}_c(\mathbf{y}) = (J_1^c(\mathbf{y}),\dots,J_{n-1}^c(\mathbf{y}))$. First, suppose that $y_1<\dots<y_n$. We set $g(\lambda):=\prod_{i=1}^n(\lambda-y_i)$. If we substitute $\mathbf{y} = (\lambda_1^+,\dots,\lambda_n^+)$ and $\alpha= c e^{-c/2}$ from \eqref{eq:delta}, we obtain 
\begin{equation}\label{eq:deltac}c^{-n}e^{nc/2}g(\lambda)+2 =: \Delta_c(\lambda),\end{equation}
where $g(\lambda) = \prod\limits_{i=1}^n (\lambda - y_i)$.  Let $\lambda_1(c)<\dots<\lambda_{2n}(c)$ be the zeros of $\Delta_c^2(\lambda)-4$; these are precisely the eigenvalues of the matrix $Q$ of \eqref{eq:Q} for $(b_1,\dots,b_n,a_1,\dots,a_n)=\varphi_{ce^{-c/2}}(c\, \mathbf{q},\mathbf{p})$. We fix $i=1\ldots, n-1$. By \eqref{eq:deltageq2}, $|\Delta_c(\mu)| \geq 2$ for any $\mu \in [\lambda_{2i}(c),\lambda_{2i+1}(c)]$. Equivalently, 
\begin{equation}
    \label{eq:g}
    |g(\mu)+2 c^ne^{-nc/2}|\ge 2 c^ne^{-nc/2} \text{ for any }\mu\in[\lambda_{2i}(c),\lambda_{2i+1}(c)].
\end{equation}
Hence, by \eqref{eq:delta}, \eqref{eq:action}, \eqref{eq:deffjc} and \eqref{eq:deltac},
\begin{equation}\label{eq:jic}
\begin{aligned}
   & J_i^c(y_1,\dots,y_n)=\frac{2}{c}\int_{\lambda_{2i}(c)}^{\lambda_{2i+1}(c)}\cosh^{-1}\frac{|c^{-n}e^{nc/2}g(\mu)+2|}{2}\,d\mu\\
    &=\frac{2}{c}\int_{\lambda_{2i}(c)}^{\lambda_{2i+1}(c)}\log\frac{ c^{-n}e^{nc/2}\left(|g(\mu)+2 c^ne^{-nc/2}|+\sqrt{g(\mu)(g(\mu)+4 c^ne^{-nc/2})}\right)}{2}\,d\mu\\
    &=(\lambda_{2i+1}(c)-\lambda_{2i}(c))(n-\frac{2n\log c}{c})
    \\&+\int_{\lambda_{2i}(c)}^{\lambda_{2i+1}(c)}\frac{2}{c}\log\left(\frac{|g(\mu)+2 c^ne^{-nc/2}|+\sqrt{g(\mu)(g(\mu)+4c^ne^{-nc/2}})}{2}\right)\,d\mu.
    \end{aligned}
\end{equation}
We set \[h_c(\mu):=\frac{2}{c}\log\left(\frac{|g(\mu)+2c^ne^{-nc/2}|+\sqrt{g(\mu)(g(\mu)+4 c^ne^{-nc/2}})}{2}\right).\]
By \eqref{eq:g}, 
\begin{equation}\label{eq:hinf}h_c(\mu)\ge -n+\frac{2n\log c}{c}.\end{equation}
Next there exist constants $M_1,M_2>0$ such that $c^ne^{-nc/2} \leq M_1$, and $|g(\mu)|\le M_2$, for all $\mu\in[\lambda_{1}(c),\lambda_{2n}(c)]$. Hence,
\begin{equation}\label{eq:hsup}
  h_c(\mu)\le \frac{2}{c}\log\left(\frac{M_1+2M_2+\sqrt{M_1(M_1+4M_2)}}{2}\right).
\end{equation}
Moreover, $h_c(\mu)\to 0$ as $c\to + \infty$ for every $\mu$. By definition of $\lambda_1(c) < \ldots < \lambda_{2n}(c)$,
\begin{equation}
    \label{eq:limp}
    \lim_{c\to + \infty}\lambda_{2i-1}(c)=\lim_{c\to +\infty}\lambda_{2i}(c)=y_i \text{  for each } i=1,\dots,n.
\end{equation}
By \eqref{eq:jic}, \eqref{eq:hinf}, \eqref{eq:hsup} and \eqref{eq:limp}, the dominated convergence theorem implies that \begin{equation}\label{eq:jlim}\lim_{c\to 
 + \infty}J_i^c(y_1,\dots,y_n)=n(y_{i+1}-y_i) \text{  for each } i=1,\dots,n.\end{equation}
By continuity of $\mathbf{J}_c$, \eqref{eq:jlim} holds for every $\mathbf{y} \in W^{n-1}_+$. Hence, so does \eqref{eq:limJ}. Moreover, since the convergence in \eqref{eq:limp} is uniform on compact sets and the bounds in \eqref{eq:hinf} and \eqref{eq:hsup} do not depend on $\mu$, by \eqref{eq:jic} the convergence in \eqref{eq:limJ} is uniform on compact sets. 
\end{enumerate}
\end{proof}

Finally, we proceed with the proof of our main result.
\begin{proof}[{\bf Proof of Theorem~\ref{thm-lag-prod-as-toric-domain}}]
We fix an open, $S_{n+1}$-symmetric, bounded and star-shaped subset $A\subset W^n$ and $\varepsilon>0$. There exists a compact set $K$ such that $(1-\varepsilon)(\mathfrak{S}^{n-1}\times_L A)\subset K\subset \mathfrak{S}^{n-1}\times_L A$. By Proposition \ref{prop:lim}, for all $c$ sufficiently large,
\[\mathbf{J}_c\circ\mathbf{F}_c((1-\varepsilon)(\mathfrak{S}^{n-1}\times_L A))\subset \mathbf{J}_c\circ\mathbf{F}_c(K)\subset \rho(A).\]
Hence, by \eqref{eq:alc}, $\Psi_c\left((1-\varepsilon)(\mathfrak{S}^{n-1}\times_L A)\right)\subset X_{\rho(A)}$ for all $c$ sufficiently large, so that 
\[(1-\varepsilon)(\mathfrak{S}^{n-1}\times_L A)\hookrightarrow X_{\rho(A)}.\]
To prove the existence of a symplectic embedding in the other direction, we observe that, by Proposition \ref{prop:lim}(b), for all $c$ sufficiently large,
\begin{equation}
    \label{eq:bounded}
   \mathbf{J}_c^{-1}(\rho(A))\subset (1+\varepsilon/2) \mathbf{J}_\infty^{-1}(\rho(A))=(1+\varepsilon/2) \left(A\cap W^{n-1}_+\right). 
\end{equation}
We set $Y_c:=\Psi^{-1}_c(X_{\rho(A)})$; by \eqref{eq:alc}, $\mathbf{F}_c(Y_c)=\mathbf{J}_c^{-1}(\rho(A))$. In particular, by \eqref{eq:bounded}, $\mathbf{F}_c(Y_c)$ bounded. Hence, by \eqref{eq:flim}, $Y_c\subset \mathfrak{S}^{n-1}\times_L W^{n-1}$ for all $c$ sufficiently large. By Proposition \ref{prop:lim}, each $\mathbf{y}\in(1+\varepsilon/2)(\bar{A}\cap W^{n-1}_+)$ has a neighborhood $U_{\mathbf{y}}$ such that \[\mathbf{F}_c^{-1}(U_{\mathbf{y}})\subset \mathfrak{S}^{n-1}\times_L (1+\varepsilon) A,\]
for all $c$ sufficiently large, where $\bar{A}$ denotes the closure of $A$. Since $\bar{A}\cap W^{n-1}_+$ is compact, there exists $c$ sufficiently large such that
\[\Psi_c^{-1}(X_{\rho(A)})=\mathbf{F}_c^{-1}\circ\mathbf{J}_c^{-1}(\rho(A))\subset \mathbf{F}_c^{-1} ((1+\varepsilon/2)(\bar{A}\cap W^{n-1}_+))\subset \mathfrak{S}^{n-1}\times_L (1+\varepsilon) A. \]
Therefore, the restriction of $\Psi^{-1}_c$ to $X_{\rho(A)}$ gives a symplectic embedding
\[X_{\rho(A)}\hookrightarrow\mathfrak{S}^{n-1}\times_L (1+\varepsilon) A\subset (1+\varepsilon)(\mathfrak{S}^{n-1}\times_L A). \]
\end{proof}

\end{document}